\documentclass[12pt]{article}
\usepackage{graphics,epsfig,graphicx,color}
\graphicspath{{../Figures/}{Figures/}}
\usepackage[caption = false]{subfig}
\usepackage{float}
\usepackage{capt-of}
\usepackage{cite}
\usepackage{enumerate}

\usepackage{amssymb,latexsym}

\usepackage{setspace}

\usepackage{amsmath}

\usepackage{mathrsfs,amsfonts}

\usepackage{amsthm,amsxtra}

\usepackage[final]{showkeys}

\usepackage{stmaryrd}

\usepackage[ruled,vlined,noresetcount]{algorithm2e}
\usepackage{algorithmicx}
\usepackage{algpseudocode}

\usepackage[openbib]{currvita}

\newif\ifPDF
\ifx\pdfoutput\undefined
	\PDFfalse
\else
	\ifnum\pdfoutput > 0
		\PDFtrue
	\else
		\PDFfalse
	\fi
\fi

\ifPDF
	\usepackage{pdftricks}
	\begin{psinputs}
		\usepackage{pstricks}
		\usepackage{pstcol}
		\usepackage{pst-plot}
		\usepackage{pst-tree}
		\usepackage{pst-eps}
		\usepackage{multido}
		\usepackage{pst-node}
		\usepackage{pst-eps}
	\end{psinputs}
\else
	\usepackage{pstricks}
\fi

\ifPDF
	\usepackage[debug,pdftex,colorlinks=true, 
	linkcolor=blue, bookmarksopen=false,
	plainpages=false,pdfpagelabels]{hyperref}
\else
	\usepackage[dvips]{hyperref}
\fi

\usepackage{fancyhdr}

\usepackage{comment}

\usepackage[toc,page]{appendix}

\usepackage{cancel}

\usepackage{esint}

\usepackage{enumitem}

\newtheorem{theorem}{Theorem}[section]

\newtheorem{remark}[theorem]{Remark}

\newcommand{\sgn}{\operatorname{sgn}}

\newcommand{\dsum}{\displaystyle\sum}

\newcommand{\dint}{\displaystyle\int}

\newcommand{\eps}{\varepsilon}

\newcommand{\bbR}{\mathbb R} \newcommand{\bbS}{\mathbb S}

 \newcommand{\bn}{\mathbf n}
 \newcommand{\bp}{\mathbf p}
\newcommand{\bq}{\mathbf q} 
 
 \newcommand{\bv}{\mathbf v} 
 \newcommand{\bx}{\mathbf x} 
\newcommand{\by}{\mathbf y} \newcommand{\bz}{\mathbf z}

 \newcommand{\cF}{\mathcal F}
\newcommand{\cG}{\mathcal G} 
 
\newcommand{\cK}{\mathcal K} \newcommand{\cL}{\mathcal L}
 \newcommand{\cN}{\mathcal N}
\newcommand{\cO}{\mathcal O}  
 \newcommand{\cR}{\mathcal R}

\newcommand{\ubar}[1]{\underline{#1}}

\newenvironment{keywords}
{\noindent{\bf Key words.}\small}{\par\vspace{1ex}}
\newenvironment{AMS}
{\noindent{\bf AMS subject classifications 2010.}\small}{\par}

\newcommand{\wt}{\widetilde}

\newcommand{\GMRES}{\textrm{GMRES}~}
\newcommand{\MINRES}{\textrm{MINRES}~}

\newcommand{\Td}{T^{\text{\tiny{dir}}}}
\newcommand{\Tf}{T^{\text{\tiny{fmm}}}}

\newcommand{\Tdg}{T^{\text{\tiny{dir}}}_{\text{\tiny{g}}}}
\newcommand{\Tfg}{T^{\text{\tiny{fmm}}}_{\text{\tiny{g}}}}
\newcommand{\Tfpg}{T^{\text{\tiny{fmm}}}_{\text{\tiny{pg}}}}
\newcommand{\Idg}{I^{\text{\tiny{dir}}}_{\text{\tiny{g}}}}
\newcommand{\Ifg}{I^{\text{\tiny{fmm}}}_{\text{\tiny{g}}}}
\newcommand{\Ifpg}{I^{\text{\tiny{fmm}}}_{\text{\tiny{pg}}}}
\newcommand{\ep}[2]{{#1}\mathrm{E}{#2}}

\title{A fast algorithm for radiative transport in isotropic media}

\author{
	Kui Ren\thanks{
		Department of Mathematics and the Oden Institute, The University of Texas, Austin, TX 78712;
		\href{mailto:ren@math.utexas.edu}{ren@math.utexas.edu}
	}
	\and
	Rongting Zhang\thanks{
		Department of Mathematics, The University of Texas, Austin, TX 78712;
		\href{mailto:rzhang@math.utexas.edu}{rzhang@math.utexas.edu}
	}
	\and
	Yimin Zhong\thanks{
		Department of Mathematics, The University of Texas, Austin, TX 78712;
		\href{mailto:yzhong@math.utexas.edu}{yzhong@math.utexas.edu}
	}
}

\begin{document}
%%%%%%%%%%%%%%%%%%%%%%%%%%%%%%%%%%%%%%%%%%%%%%%%%%%%%%%%%%%%%%%%%%
%%%%%%%BEGIN DOCUMENT %%%%BEGIN DOCUMENT %%%%BEGIN DOCUMENT%%%%%%%
%%%%%%%%%%%%%%%%%%%%%%%%%%%%%%%%%%%%%%%%%%%%%%%%%%%%%%%%%%%%%%%%%%

\maketitle

%\tableofcontents

%%%%%%%%%%%%%%%%%%%%%%%%%%%
%%%%%%%%ABSTRACT%%%%%%%%%%%
%%%%%%%%%%%%%%%%%%%%%%%%%%%
\begin{abstract}
Constructing efficient numerical solution methods for the equation of radiative transfer (ERT) remains as a challenging task in scientific computing despite of the tremendous development on the subject in recent years. We present in this work a simple fast computational algorithm for solving the ERT in isotropic media. The algorithm we developed has two steps. In the first step, we solve a volume integral equation for the angularly-averaged ERT solution using iterative schemes such as the GMRES method. The computation in this step is accelerated with a fast multipole method (FMM). In the second step, we solve a scattering-free transport equation to recover the angular dependence of the ERT solution. The algorithm does not require the underlying medium be homogeneous. We present numerical simulations under various scenarios to demonstrate the performance of the proposed numerical algorithm for both homogeneous and heterogeneous media.
\end{abstract}

%%%%%%%%%%%%%%%%%%%%%%%%%%%
%%%%%%%%%KEYWORDS%%%%%%%%%%
%%%%%%%%%%%%%%%%%%%%%%%%%%%


%%%%%%%%%%%%%%%%%%%%%%%%%%
%%%   AMS or PACS   %%%%%%
%%%%%%%%%%%%%%%%%%%%%%%%%%


%%%%%%%%%%%%%%%%%%%%%%%%%%%%%%%%%%%%%%%%%%%%%%%%%%%%%%%%%%%%%%%%%%
%%%%%% BEGINNING TEXT %%% BEGINNING TEXT %%% BEGINNING TEXT %%%%%%
%%%%%%%%%%%%%%%%%%%%%%%%%%%%%%%%%%%%%%%%%%%%%%%%%%%%%%%%%%%%%%%%%%

%\begin{comment}

%1. quadrature rules

%2. analytical solutions

%3. preconditioning

%4. specialty of the kernel (different from screened Colomb)

%\end{comment}

%%%%%%%%%%%%%%%%%%%%%%%%%%%%%%%%%%%%%%%%%%%%%%%%%%%%%%%%%%%%%%%%%%
%%%%%%%%%%%%%%%%%%%%%%%%%%%%%%%%%%%%%%%%%%%%%%%%%%%%%%%%%%%%%%%%%%
\section{Introduction}
\label{SEC:Intro}
%%%%%%%%%%%%%%%%%%%%%%%%%%%%%%%%%%%%%%%%%%%%%%%%%%%%%%%%%%%%%%%%%%
%%%%%%%%%%%%%%%%%%%%%%%%%%%%%%%%%%%%%%%%%%%%%%%%%%%%%%%%%%%%%%%%%%

This work is concerned with the numerical solution of the steady-state equation of radiative transfer (ERT) with spatially dependent physical coefficients and \emph{isotropic} scattering kernel~\cite{DaLi-Book93-6,Graziani-Book06,LeMi-Book93}:
\begin{equation}\label{EQ:ERT}
\begin{array}{rcll}
\bv \cdot \nabla \Phi(\bx, \bv) + \mu(\bx) \Phi(\bx, \bv) - \mu_s(\bx) \dint_{\bbS^{d-1}} \Phi(\bx, \bv') d\bv' &=& f(\bx),& \mbox{in}\ \Omega\times\bbS^{d-1}\\
\Phi(\bx, \bv) &=& 0,& \mbox{on}\ \Gamma_{-}
\end{array}
\end{equation}
where $\Omega\subseteq\bbR^d$ ($d=2,3$) is a bounded domain with smooth boundary $\partial\Omega$, $\bbS^{d-1}$ is the unit sphere in $\bbR^d$, and $\Gamma_{-}=\{(\bx,\bv):\ (\bx, \bv)\in \partial\Omega\times\bbS^{d-1}\ \mbox{s.t.}\ \bn(\bx)\cdot \bv<0 \}$ ($\bn(\bx)$ being the unit outer normal vector at $\bx\in\partial\Omega$) is the incoming part of the phase space boundary. For the only reason of simplifying the presentation, we have assumed that there is no incoming source on the boundary. Moreover, we have assumed that the internal source $f$ is only a function of the spatial variable. In fact, this is not needed either for our algorithm to work; see more discussions in Section~\ref{SEC:Concl}.

The equation of radiative transfer is a popular model for describing the propagation of particles in complex media. It appears in many fields of science and technology, ranging from classical fields such as nuclear engineering~\cite{Larsen-NSE88,Larsen-TTSP88,Mokhtar-Book97}, astrophysics~\cite{HeGr-AJ41,CeBaBeAi-TTSP99,WaUe-ASS89}, and remote sensing~\cite{BaRe-IP05,SpKuCh-JQSRT01}, to modern applications such as biomedical optics~\cite{Arridge-IP99,DiRe-JCP14,KiMo-IP06,Ren-CiCP10,ReBaHi-SIAM06,ReAbBaHi-OL04}, radiation therapy and treatment planning~\cite{Boman-Thesis07,JiScPaJi-PMB12,TiHeToSiAlPyUl-PMB08}, and imaging in random media~\cite{BaRe-SIAM08,BaCaLiRe-IP07,BoGa-IP16,RyPaKe-WM96}. The coefficients $\mu(\bx)$ and $\mu_s(\bx)$ have different physical meanings in different applications. In general, the coefficient $\mu_s(\bx)$ measures the strength of the scattering of the underlying medium at $\bx$, while $\mu_a(\bx)\equiv \mu(\bx)-\mu_s(\bx)$ measures the strength of the physical absorption of the medium. The coefficient $\mu(\bx)$ measures the total absorption at $\bx$ due to both the physical absorption and absorption caused by scattering, that is the loss of particles from the current traveling direction into other directions due to scattering.

Numerical methods for solving the equation of radiative transfer has been extensively studied, see for instance~\cite{DaLi-Book93-6,Graziani-Book06,LaPaSe-Book03,LeMi-Book93,RoViBo-PNE11} and references therein for an overview. Besides Monte Carlo type of methods that are based on stochastic representation of the ERT~\cite{BhSp-JCP07,DeThUr-JCP12,Edstrom-SIAM05,HaSpVe-SIAM07,UeLa-JCP98}, many different deterministic discretization schemes have been proposed~\cite{AdNo-JCP98,AnLa-JCP01,Asadzadeh-SIAM98,BoLaAd-JCP92,DeVo-JCP02,DuKl-JCP02,FrKlLaYa-JCP07,GaZh-TTSP09,GrSc-JCP11,GuKa-SIAM10,Hermeline-JCP16,JiPaTo-SIAM00,KaRa-SIAM14,KiBeGoMo-JCP10,KiMo-SIAM02,LaMcHa-JCP16,LaThKlSeGo-JCP02,Lesaint-FEP86,MaReSt-SIAM00,MoRaAdKa-TTSP13,RaGuKa-JCP12,TuFrDuKl-JCP04} and numerous iterative schemes, as well as preconditioning strategies, have been developed to solve the discretized systems; see for instance~\cite{AdLa-PNE02,BrMoRa-JCP14,GoLi-JCP12,OlDe-PNE98,PaHo-ANE02} and references therein.

There are many challenging issues in the numerical solutions of the equation of radiative transfer. One of such challenges is the high-dimensionality involved. The ERT is posed in phase space, meaning that the main unknown in the equation, in steady state, depends on both the spatial variable $\bx\in\Omega$ and the angular variable $\bv\in\bbS^{d-1}$. In the spatial three-dimensional case, the unknown $\Phi$ depends on five variables, three in the spatial domain and two in the angular domain. This poses significant challenges in terms of both solution speed and storage.

In this work, we propose a new method to solve the ERT in isotropic media, that is, media whose physical coefficients and the scattering kernel do not depend on the angular variable $\bv$, i.e., the media absorb and scatter particles in the same manner for all directions. Our method is based on the observation that when the underlying medium is isotropic, the angularly averaged ERT solution, $\int_{\bbS^{d-1}} \Phi(\bx,\bv)\,d\bv$, satisfies a Fredholm integral equation of the second type. This integral equation can be solved, using a fast multiple method, for $\int_{\bbS^{d-1}} \Phi(\bx,\bv)\,d\bv$. Once this is done, we can plug $\int_{\bbS^{d-1}} \Phi(\bx,\bv)\,d\bv$ into the ERT~\eqref{EQ:ERT} to solve for $\Phi$ itself.

The rest of this paper is organized as follows. In Section~\ref{SEC:Integral}, we re-formulate the ERT~\eqref{EQ:ERT} into a Fredholm integral equation of the second type for the unknown $\int_{\bbS^{d-1}} \Phi(\bx,\bv)\,d\bv$. We then propose in Section~\ref{SEC:FMM} a numerical procedure for solving the ERT based on this integral formulation and implement an interpolation-based fast multipole method~\cite{FoDa-JCP09} to solve the integral equation. Important issues on the implementation of our method are discussed in Section~\ref{SEC:Impl}. In Section~\ref{SEC:Num} we present some numerical tests for the algorithm that we developed. Concluding remarks are then offered in Section~\ref{SEC:Concl}.

%%%%%%%%%%%%%%%%%%%%%%%%%%%%%%%%%%%%%%%%%%%%%%%%%%%%%%%%%%%%%%%%%%
%%%%%%%%%%%%%%%%%%%%%%%%%%%%%%%%%%%%%%%%%%%%%%%%%%%%%%%%%%%%%%%%%%
\section{Integral formulation}
\label{SEC:Integral}
%%%%%%%%%%%%%%%%%%%%%%%%%%%%%%%%%%%%%%%%%%%%%%%%%%%%%%%%%%%%%%%%%%
%%%%%%%%%%%%%%%%%%%%%%%%%%%%%%%%%%%%%%%%%%%%%%%%%%%%%%%%%%%%%%%%%%

Our algorithm is based on the integral formulation of the ERT~\eqref{EQ:ERT}. This is a well-developed subject. We refer to~\cite{DaLi-Book93-6} for more details. To present the formulation, let us first introduce a function $q(\bx)$ defined as 
\begin{equation*}
	q(\bx) :=  \mu_s(\bx) \dint_{\bbS^{d-1}} \Phi(\bx, \bv') d\bv' + f(\bx). 
\end{equation*}
We can then rewrite the equation of radiative transfer, using the method of characteristics, into the following integral form~\cite{DaLi-Book93-6}:
\begin{equation}\label{EQ:ERT Intgrl}
	\Phi(\bx, \bv) = \int_0^{\tau(\bx, \bv)} \exp\left(-\int_0^\ell \mu(\bx-\ell'\bv )d\ell'\right) q(\bx-\ell \bv) d\ell .
\end{equation}
Here $\tau(\bx,\bv)$ is the distance it takes for a particle to go from $\bx$ to reach the domain boundary $\partial\Omega$ in the $-\bv$ direction:
\[
	\tau(\bx,\bv)=\sup\{\ell:\ \bx-\ell'\bv\in\Omega\ \mbox{for}\ 0\le \ell'< \ell\}.
\]
The integral formulation in~\eqref{EQ:ERT Intgrl} is classical and has been used to derive many theoretical results and numerical methods on the ERT~\cite{DaLi-Book93-6,Mokhtar-Book97}. 

The most crucial step of our algorithm is to integrate the integral formulation~\eqref{EQ:ERT Intgrl} again over $\bbS^{d-1}$ to obtain an integral equation for the local density $U(\bx)$ defined by
\begin{equation*}
	U(\bx) :=  \dint_{\bbS^{d-1}} \Phi(\bx, \bv) d\bv.
\end{equation*}
The result is a Fredholm integral equation of the second type. It reads
\begin{equation}\label{EQ:ERT Intgrl U}
	U(\bx) = KU(\bx) + K(\mu_s^{-1}f)(\bx), 
\end{equation}
where the linear integral operator $K$ is defined as
\begin{equation*}
	K g(\bx) =\dint_{\bbS^{d-1}} \int_0^{\tau(\bx, \bv)} \mu_s(\bx-\ell\bv) \exp\left(-\int_0^\ell \mu(\bx - \ell'\bv) d\ell' \right)  g(\bx- \ell\bv) d\ell d\bv .
\end{equation*}
To simplify the expression for $K$, let $\by = \bx - \ell \bv$, and define the function $E(\bx, \by)$
\begin{equation*}
	E(\bx, \by) = \exp\left(-\int_0^{|\by-\bx|} \mu (\bx -  \ell' \dfrac{\bx-\by}{|\bx-\by|}) d\ell' \right),
\end{equation*}
which is nothing but the total absorption along the line segment between $\bx$ and $\by$. We can then express the integral operator $K$ as
\begin{equation}\label{EQ:Operator K}
	Kg(\bx) = \int_\Omega \cK(\bx,\by) g(\by) d\by
\end{equation}
where the integral kernel $\cK$ is defined as
\begin{equation}\label{EQ:Kernel}
	\cK(\bx, \by) = \frac{1}{|\bbS^{d-1}|} \frac{\mu_s(\by) E(\bx,\by)}{|\bx -\by|^{d-1}}
\end{equation}
with $|\bbS^{d-1}|$ the surface area of the unit sphere $\bbS^{d-1}$. $|\bbS^{d-1}|=2\pi$ when $d=2$ and $|\bbS^{d-1}|=4\pi$ when $d=3$. In the case where $\mu$ and $\mu_s$ are independent of the spatial variable, the integral kernel $\cK$ simplifies to
\begin{equation}\label{EQ:Kernel Cons}
	\cK(\bx,\by) = \frac{1}{|\bbS^{d-1}|} \frac{\mu_s e^{-\mu |\bx-\by|}}{|\bx - \by|^{d-1}}.
\end{equation}
This integral kernel has the same form as the Yukawa potential (also called a screened Coulomb potential) when $d=2$~\cite{GrHu-JCP02}. 

The algorithm we propose here is based on the integral formulation of the ERT for the variable $U$ that we derived in~\eqref{EQ:ERT Intgrl U}. The integral operator $K$ is compact since the kernel function $\cK(\bx, \by)$ is weakly singular.
From (\ref{EQ:ERT Intgrl U}), we can obtain that
\begin{equation}\label{EQ:ERT Intgrl FMM}
	(I - K)U(\bx) = \phi(\bx),
\end{equation}
where $\phi(\bx)\equiv K(\mu_s^{-1} f)(\bx)$. The operator $(I -K)$ is a Fredholm operator, and by Fredholm alternative theorem and the fact that the ERT admits only the zero solution when $f\equiv 0$, see for instance~\cite{DaLi-Book93-6}, we conclude that there is a unique solution to~\eqref{EQ:ERT Intgrl FMM}.

Let us finish this section with the following important observation. The kernel~\eqref{EQ:Kernel} for the volume integral equation that we derived here takes the same form in the cases of homogeneous (i.e. $\mu$ and $\mu_s$ do not depend on spatial variable) and inhomogeneous (i.e. $\mu$ and $\mu_s$ depend on spatial variable) media. This means that the algorithm that we present in the next sections work for both homogeneous and inhomogeneous media. In the case of homogeneous media or inhomogeneous media that is analytically known, the evaluation of the kernel $\cK(\bx, \by)$ can be analytically performed. When the coefficient $\mu$ is only given on a collection of points in the domain, the kernel needs to be evaluated with a numerical quadrature rule for $E(\bx, \by)$. The computation of $\cK(\bx, \by)$ in this case is therefore more expensive.

%%%%%%%%%%%%%%%%%%%%%%%%%%%%%%%%%%%%%%%%%%%%%%%%%%%%%%%%%%%%%%%%%%
%%%%%%%%%%%%%%%%%%%%%%%%%%%%%%%%%%%%%%%%%%%%%%%%%%%%%%%%%%%%%%%%%%
\section{A fast multipole based algorithm}
\label{SEC:FMM}
%%%%%%%%%%%%%%%%%%%%%%%%%%%%%%%%%%%%%%%%%%%%%%%%%%%%%%%%%%%%%%%%%%
%%%%%%%%%%%%%%%%%%%%%%%%%%%%%%%%%%%%%%%%%%%%%%%%%%%%%%%%%%%%%%%%%%

Our strategy of solving the ERT~\eqref{EQ:ERT} is to first solve for $U$ and then solve for $\Phi$ from $U$. The main solution procedure is as follows.

\begin{algorithm}[H]
    \SetAlgoLined
    \KwData{coefficients $\mu(\bx)$, $\mu_s(\bx)$, source $f(\bx)$}
    \KwResult{Solution $\Phi(\bx, \bv)$ to ERT~\eqref{EQ:ERT} }
{\bf{S.1}~}evaluate the source function $\phi(\bx)\equiv K(\mu_s^{-1}f)(\bx)$ analytically, or by:\
\begin{enumerate}[label ={(\roman*)}]
    \item solving the following scattering-free transport equation for $u$:
    \begin{equation*}\label{EQ:ERT Free1}
    \begin{array}{rcll}
    \bv \cdot \nabla u(\bx, \bv) + \mu(\bx) u(\bx, \bv) &=& \mu_s^{-1}(\bx)f(\bx),& \mbox{in}\ \Omega\times\bbS^{d-1}\\
    u(\bx, \bv) &=& 0,& \mbox{on}\ \Gamma_{-}
    \end{array}
    \end{equation*}
    \item evaluating $\phi(\bx)=\int_{\bbS^{d-1}} u(\bx,\bv) d\bv$.
\end{enumerate}
 {\bf{S.2}~}   Use a Krylov subspace method, such as the \GMRES algorithm~\cite{Saad-Book03}, to solve the integral equation~\eqref{EQ:ERT Intgrl FMM} for $U$.
 
 \vspace{0.2cm}
 {\bf{S.3}~}  Recover the ERT solution $\Phi$ by 
\begin{enumerate}[label ={(\roman*)}]
    \item evaluating the source $Q(\bx)=\mu_s(\bx) U(\bx)+f(\bx)$;
    \item solving the following scattering-free transport equation for $\Phi$:
    \begin{equation*}\label{EQ:ERT Free}
    \begin{array}{rcll}
    \bv \cdot \nabla \Phi(\bx, \bv) + \mu(\bx) \Phi(\bx, \bv) &=& Q(\bx),& \mbox{in}\ \Omega\times\bbS^{d-1}\\
    \Phi(\bx, \bv) &=& 0,& \mbox{on}\ \Gamma_{-} 
    \end{array}
    \end{equation*}
\end{enumerate}
\caption{General Solution Procedure}
\label{ALG:FMM}
\end{algorithm}

The solution of the scattering-free transport equations in the first and last steps can be done efficiently with a fast sweeping method such as that in~\cite{GaZh-TTSP09} or even analytically in special cases. The solution of the integral equation in the second step is more or less straightforward since the integral kernel is only weakly singular. Nevertheless, in the rest of the paper, we still present some numerical evidences to demonstrate the performance of our algorithm.

Let us remark that one feature of the above method for solving the ERT~\eqref{EQ:ERT} is that it \emph{does not} require an explicit discretization over the angular variable. It is clear that the main computational cost of the algorithm is on the solution of the integral equation~\eqref{EQ:ERT Intgrl FMM} which involves only the spatial variable. Therefore, besides the solution of the scattering-free transport equation, the computational complexity of the algorithm \emph{does not} scale with the size of the angular discretization. In many applications, the main quantities of interests is the local density $U(\bx)$, not $\Phi(\bx,\bv)$. In these cases, the ${\bf S.3}$ step of Algorithm~\ref{ALG:FMM} is not necessary, and the computational complexity of the algorithm therefore is completely independent of the angular discretization beside the one transport sweep in the construction of $\phi$ for ~\eqref{EQ:ERT Intgrl FMM}. For the same reason, the storage requirement of the algorithm also depends only on the spatial discretization.

%%%%%%%%%%%%%%%%%%%%%%%%%%%%%%%%%%%%%%%%%%%%%%%%%%%%%%%%%%%%%%%%%%
\subsection{Discretization}
\label{SEC:DISC}
%%%%%%%%%%%%%%%%%%%%%%%%%%%%%%%%%%%%%%%%%%%%%%%%%%%%%%%%%%%%%%%%%%

There are many existing methods of the discretization for integral equations with weakly singular kernels, see for instance~\cite{Hochstadt-Book89,Kress-Book99,Vainikko-Book93, BrGi-JCP12} and references therein. It has been shown in~\cite{Vainikko-Book93} that when the kernel $\cK(\bx, \by)$ is differentiable on $\Omega\times\Omega \backslash \{ \bx = \by \}$, then the solution to volume integral equation~\eqref{EQ:ERT Intgrl FMM} \emph{only} belongs to H\"{o}lder continuous space, $U\in C^{0,\alpha}(\overline{\Omega})$, $\forall \alpha\in (0,1)$.
In fact, it is easy to verify that 
$$|U(\bx) - U(\by)|\le \cO\left( |\bx - \by|\log (|\bx - \by|)\right).$$
This deficiency in regularity comes from the boundary effect. The vacuum boundary condition in~\eqref{EQ:ERT} implies that $\mu_s \equiv 0$ outside $\Omega$, which imposes a jump across the boundary $\partial\Omega$. 
Therefore high order discretization schemes for integral equations like~\eqref{EQ:ERT Intgrl FMM} usually require special treatments to the boundary effect~\cite{Vainikko-Book93} unless $\Omega$ does not have a boundary.

For simplicity, in the following we consider the \emph{piecewise constant collocation method} (PCCM)~\cite{Vainikko-Book93} with the following assumptions: (i) the domain $\Omega$ is convex with Lipschitz boundary $\partial\Omega$; (ii) the coefficients $\mu, \mu_s\in C^2(\Omega)$; and (iii) the source function $f\in C^2(\Omega)$. The discretization is constructed as follows:
\begin{enumerate}
    \item Partition of $\Omega$. For a small $h > 0$, we partition the spatial domain $\Omega$ into two parts: $\Omega_b^h$ and $\Omega_{i}^h$, where
    \begin{equation}
    \Omega_{b}^h := \{\bx\in\Omega : \text{dist}(\partial\Omega, \bx) \le h^2 \}\;\text{ and }\; \Omega_i^h = \Omega \backslash \Omega_b^h, 
    \end{equation}
     Take a discretization $\{T_{j,h}\}_{j=1}^N$ of $\Omega$, that is, $T_{j,h}\cap T_{j',h}=\emptyset$ $\forall j\neq j'$ and $\Omega = \bigcup_{j=1}^{N} T_{j,h}$, such that (a) $\text{diam}(T_{j,h}) \le h$ $\forall j$ and (b) $T_{j,h}\cap \Omega_i^h\neq \emptyset$ $\forall j$ (which is only saying that there is no cell $T_{j,h}$ that is completely inside $\Omega_b^h$, an assumption that can be easily satisfied since the thickness of $\Omega_i^h$ is of order $h^2$). It is then clear that $N \simeq \cO(h^{-d})$. For any $1\le j\le N$, if $T_{j,h}\cap \Omega_b^h\neq \emptyset$, we set $T'_{j,h} := T_{j,h} \cap \Omega_i^{h}$ when it is not empty.
     \item Collocation Points. For each cell in the discretization, we locate the collocation point $\bx_j \in T_{j,h}$ by
     \begin{enumerate}
         \item If $T_{j,h}\subset\Omega_i^h$, then we choose $\bx_j$ as the centroid point
         \begin{equation}
         \bx_j = \frac{1}{|T_{j,h}|} \int_{T_{j,h}} \bz d\bz.
         \end{equation}
         \item If $T_{j,h}\cap \Omega_{b}^h\neq \emptyset$, then we choose arbitrary $\bx_j\in T'_{j,h}$.
        \end{enumerate}
\end{enumerate}
The simplest example of the above discretization is to use a uniform grid $\cG$ with cell size of $h$. For a cell $T_{j,h}\subset \cG$ contained in $\Omega$, we choose its centroid point as the collocation point. For a boundary-incident cell $T_{j,h}\subset \cG$ such that $T_{j,h}\cap \partial\Omega\neq \emptyset$, we replace the cell $T_{j,h}$ by the intersection $T_{j,h}' = T_{j,h}\cap\Omega$ and choose an arbitrary point in $T'_{j,h}$ as the collocation point. When the boundary $\partial\Omega$ is $C^2$, the boundary part $\partial\Omega\cap T_{j,h}$ can be approximated using a tangent plane or secant plane. The omitted measure in this case is at the order of $\cO(h^2)$. 

%%%%%%%%%%%%%%%%%%%%%%%%%%%%%%%%%%%%%%%%%%%%%%%%%%%%%%%%%%%%%%%%%%
\subsection{Linear system from discretization}
\label{SEC:LINEAR SYS}
%%%%%%%%%%%%%%%%%%%%%%%%%%%%%%%%%%%%%%%%%%%%%%%%%%%%%%%%%%%%%%%%%%

We then represent the piecewise constant solution $\bar{U}(\bx)$ by
\begin{equation}
\bar{U}(\bx) = \sum_{j=1}^{N} u_{j,h} \chi_{j,h}(\bx),\quad \chi_{j,h}(\bx) = \begin{cases}
1,\quad \bx \in T_{j,h} \\
0,\quad \bx \notin T_{j,h}
\end{cases}
\end{equation}
Replacing $U$ by $\bar{U}$ in the integral equation~\eqref{EQ:ERT Intgrl FMM}, we obtain the discretized linear equation for $\bar{U}(\bx_j) = u_{j,h}$:
\begin{equation}\label{EQ: DISCRETIZED}
\bar{U}(\bx_{j}) = \sum_{k=1}^{N} K_{jk}\omega_k \bar{U}(\bx_{k})  + \phi(\bx_j),\ \ j=1,\dots, N
\end{equation}
where 
\begin{equation}\label{EQ: KIJ}
K_{jk} = \frac{1}{|T_{k,h}|}\int_{T_{k, h}} \cK(\bx_{j}, \by) d\by,\quad \omega_k = |T_{k,h}|.
\end{equation}
We consider three mostly used approaches in the evaluation of the elements of the $K$ matrix in the above linear system~\eqref{EQ: DISCRETIZED}:
\begin{enumerate}[leftmargin = 1.8cm]
    \item[(M-i)] Integrate $\int_{T_{k, h}} \cK(\bx_j, \by) d\by$ in~\eqref{EQ: KIJ} analytically. This is hard to implement for general discretizations. When the cells $T_{k,h}$ are of regular shapes, for instance cubes or simplexes, we can obtain closed form evaluations.
    \item[(M-ii)] Rewrite $\cK(\bx_j, \bx_k) = a(\bx_j, \bx_k) / |\bx_j - \bx_k|^{d-1}$ with $a(\cdot, \cdot)$ twice differentiable on $\Omega\times\Omega \backslash \{ \bx = \by\}$. We then make the approximation $\int_{T_{k, h}} \cK(\bx_j, \by) d\by\approx a(\bx_j, \bx_k)\int_{T_{k,h}}|\bx_j - \by|^{1-d} d\by$. It is natural to take $a(\bx, \by) = |\bbS^{d-1}|^{-1}E(\bx, \by)\mu_s(\by)$ for our kernel in~\eqref{EQ:Kernel}. 
    \item[(M-iii)]Use the approximation $\int_{T_{k, h}} \cK(\bx_j, \by) d\by \simeq |T_{k,h}| \cK(\bx_j, \bx_k)$ when $j\neq k$. For the singular integral at $j=k$, one can compute it explicitly or simply ignore it.
\end{enumerate}
It is obviously that the approach (M-iii) has the lowest accuracy among the three approaches. However, it is the easiest approach to implement in practice. To be more precise, it can be shown following~\cite[Theorem 5.1]{Vainikko-Book93} that the numerical errors for the above discretization schemes are respectively:
\begin{equation}\label{EQ: ERROR}
\begin{aligned}
\max_{1\le j \le N} |U(\bx_{j}) - \bar{U}(\bx_{j})| &\le \cO(h^2(1+\log |h|))\quad &&\text{if using (M-i)},\\
\max_{1\le j \le N} |U(\bx_{j}) - \bar{U}(\bx_{j})| &\le \cO(h^2(1+\log |h|))\quad &&\text{if using (M-ii)},\\
\max_{1\le j \le N} |U(\bx_{j}) - \bar{U}(\bx_{j})| &\le \cO(h)\quad  &&\text{if using (M-iii)},
\end{aligned}
\end{equation}
where the constants in the estimates would depend on the coefficients $\mu$, $\mu_s$ and the source function $\phi$ which are all assumed to be smooth enough, at least in the class of $C^2(\Omega)$.

With a slight abuse of notation, we write the linear system~\eqref{EQ: DISCRETIZED} again in the form
\begin{equation}\label{EQ: LINEAR SYSTEM}
(I-K)\bar{U}=\phi,
\end{equation}
where the integral kernel matrix $K = [K_{jk}\omega_k]_{N\ge j,k \ge 1}$,  the vectors  $\bar U= [\bar U(\bx_j)]_{N\ge j \ge 1}$ and $\phi = [\phi(\bx_j)]_{N\ge j \ge 1}$.

\begin{remark}[Evaluating Elements of $K$]\label{REM: NEW KERNEL}
     In our numerical implementation, we take uniform discretizations where the cells $\{ T_{j,h} \}_{1\le j\le N}$ are identical, e.g. hypercubes or hyperrectangles. This simplifies the evaluation of the elements of the $K$ matrix for approaches (M-i) and (M-ii). For instance, in approach (M-ii), after the approximation of the kernel, we have
     \begin{equation}
     K_{jk} \approx \dfrac{a(\bx_j, \bx_k)}{|T(\bx_k)|} \int_{T(\bx_k)} |\bx_j - \bz|^{1-d} d\bz
     \end{equation}
     where $T(\bx_k)$ is a region identical to all the cells with $\bx_k$ as its centroid. The computation of $\int_{T(\bx_k)} |\bx_j - \bz|^{1-d} d\bz$ can be done with Fourier transform analytically or numerically. Consider the two-dimensional case ($d=2$), and let $T_{j,h}$ be identical and square. Let $T(\by)$ be a square centered at $\by= (y_1, y_2)$ with side length of $h$. Let $\bx = (x_1, x_2)$, $t_1 = y_1 - x_1$, and $t_2 = y_2 - x_2$. It is then easy to verify that
     \begin{equation}
     \int_{T(\by)} |\bx - \bz|^{1-d} d\bz =  \sum_{i=-1}^1\sum_{j=-1}^1 ij F\left(t_1 + i\frac{h}{2} ,\; t_2 + j\frac{h}{2}\right) 
     \end{equation} 
     with the function $F(r, s)$ given by
     \begin{multline}
     F(r, s) = \sgn(r) \sgn(s) \Big( |r|\log(|s| + \sqrt{r^2 + s^2})  \\ + |s|\log(|r| +\sqrt{r^2 + s^2})-|r|\log |r| -|s|\log |s| \Big).
     \end{multline} 
     This calculation works for any $(\bx, \bz)$ pair over $\Omega\times\Omega$ and any $h > 0$.
     
    In the same spirit, if we redefine a kernel $\wt\cK(\bx, \by)$ as
     \begin{equation}
     	\wt \cK(\bx, \by) = a(\bx, \by) \int_{T(\by)} |\bx - \bz|^{1-d} d\bz
     \end{equation}
     and replace the kernel $\cK$ in (M-iii) with the new kernel $\wt\cK$ (which does not have singularity at $\bx=\by$ since neither the function $a(., .)$ nor the function $F(., .)$ has), we achieve a better accuracy, i.e. $\cO(h^2\log|h|)$, for the corresponding discretization.
\end{remark}

%%%%%%%%%%%%%%%%%%%%%%%%%%%%%%%%%%%%%%%%%%%%%%%%%%%%%%%%%%%%%%%%%%
\subsection{Fast multipole method}
\label{SEC: FMM}
%%%%%%%%%%%%%%%%%%%%%%%%%%%%%%%%%%%%%%%%%%%%%%%%%%%%%%%%%%%%%%%%%%

To solve the linear system $(I-K)\bar{U}=\phi$ with a \GMRES or \MINRES algorithm, we need to evaluate the matrix-vector product of the form $(I-K)\bar{U}$ for different vectors $\bar{U}$. Therefore, the main computational cost will be dominated by the cost of the evaluations of $K\bar{U}$. Direct evaluation of such a summation takes $\cO(N^2)$ operations in general. In this work, we use the fast multipole method (FMM), originally developed by Greengard and Rokhlin~\cite{GrRo-JCP87}, to accelerate the evaluation of this matrix-vector product. For the simplicity of implementation, we use an interpolation-based FMM that was proposed by Fong and Darve in~\cite{FoDa-JCP09}. Other efficient implementations of FMM, see for instance~\cite{ChGrRo-JCP99,ChJiMi-Book01,GrRo-AN97,MaRo-SIAM07,YiBiZo-JCP03} and references therein, may also be applied to our problem here. This will be a future work.

In our implementation of the Fong-Darve FMM algorithm~\cite{FoDa-JCP09}, we follow the standard multilevel approach with k-d tree structure. The idea in~\cite{FoDa-JCP09} is based on the Chebyshev interpolation for the far-field interactions. 
Let $T_k(x)$ be the first-kind Chebyshev polynomial of degree $k$ defined on $[-1, 1]$. Define the interpolation function
\begin{equation}\label{EQ: CHEB}
S_n(\bp,\bq) = \prod_{i = 1}^d\left(\frac{1}{n} + \frac{2}{n}\sum_{k=1}^{n-1}T_k(p_i)T_k(q_i)\right)
\end{equation}
with the conventions $\bp=(p_1,\cdots,p_d)\in [-1,1]^d$ and $\bq=(q_1,\cdots,q_d)\in[-1,1]^d$.
Take the hyperrectangles $X_1 = \prod_{i=1}^d [\ubar{a}_i, \bar{a}_i]\subset\Omega$ and $X_2 = \prod_{i=1}^d [\ubar{b}_i, \bar{b}_i]\subset \Omega$, and assume that $X_2$ stays in the far-field of $X_1$. We define the linear transformations $\cL: [-1, 1]^d \mapsto X_1$ and $\cR: [-1,1]^d \mapsto X_2$ as the following
\begin{equation}
\begin{aligned}
\cL \bp &= \left( \frac{\ubar{a}_1 + \bar{a}_1}{2} + \frac{\bar{a}_1 - \ubar{a}_1}{2}p_1,\dots, \frac{\ubar{a}_d + \bar{a}_d}{2} + \frac{\bar{a}_d - \ubar{a}_d}{2}p_d\right),\\ \cR \bq &= \left( \frac{\ubar{b}_1 + \bar{b}_1}{2} + \frac{\bar{b}_1 - \ubar{b}_1}{2}q_1,\dots, \frac{\ubar{b}_d +\bar{b}_d}{2} + \frac{\bar{b}_d - \ubar{b}_d}{2}q_d\right).
\end{aligned}
\end{equation}
These linear transforms map the standard Cheyshev points in $[-1,1]^{d}$ to the \emph{scaled} Chebyshev points in $X_1$ and $X_2$ respectively.
If the two-variable kernel $\cK(\cdot, \cdot)$ is \emph{smooth} enough on $X_1 \times X_2$, then it can be well approximated by the following interpolation formula~\cite{DuGuRo-SIAM96,FoDa-JCP09}:
\begin{equation}\label{EQ:Kernel Approx}
\begin{aligned}
	\cK(\bx_i,\bx_j) \approx \dsum_{m=1}^{n^d} \dsum_{m'=1}^{n^d} S_n( \cL^{-1} \bx_i, \wt\bp_m) \cK( \cL\wt\bp_m, \cR \wt\bq_{m'})  S_n( \cR^{-1} \bx_j, \wt\bq_{m'}),\quad (\bx_i, \bx_j)\in X_1\times X_2
\end{aligned}
\end{equation}
where $\wt\bp_m, \wt \bq_{m'}\in Z:= \{\wt\bz_{k}\}_{k=1}^{n^d}\subset [-1,1]^d$, $Z$ being the set of $d$-dimensional Chebyshev interpolation points which are computed by the $d$-dimensional tensor-product of the $n$-th order Chebyshev points on $[-1, 1]$. Then far-field contribution between $X_1$ and $X_2$ for matrix-vector product $K\bar{U}$, for the approach (M-iii) in Section~\ref{SEC:LINEAR SYS}, can be represented by
\begin{equation}\label{EQ: SUMMATION}
\varphi(\bx_i) = \sum_{m=1}^{n^d} S_n(\cL^{-1} \bx_i, \wt\bp_m) \sum_{m'=1}^{n^d} \cK(\cL\wt{\bp}_m, \cR\wt \bq_{m'} ) \sum_{\bx_j \in X_2} S_n(\cR^{-1} \by, \wt\bq_{m'}) \bar{U}(\bx_j) \omega_j,\quad \forall \bx_i\in X_1.
\end{equation}
The above summation consists of 3 steps: 
\begin{enumerate}
    \item Evaluation $f_1(\wt\bq_{m'}) = \sum_{\bx_j \in X_2} S_n(\cR^{-1} \by, \wt\bq_{m'}) \bar{U}(\bx_j) \omega_j$ for all Chebyshev points $\wt\bq_{m'}$ in $[-1,1]^d$. 
    \item Compute the interactions between Chebyshev points $\wt\bp_m$ and $\wt\bq_{m'}$ through 
    \begin{equation}\label{EQ: M2L}
     f_2(\wt\bp_{m}) = \sum_{m'=1}^{n^d} \cK(\cL\wt{\bp}_m, \cR\wt \bq_{m'} ) f_1(\wt\bq_{m'})  
    \end{equation}
    for all Chebyshev points $\wt\bp_m\in [-1,1]^d$.
    \item Interpolate back to the collocation points $\bx_i \in X_1$ with
\begin{equation}\label{EQ: L2L}
\varphi(\bx_i) = \sum_{m=1}^{n^d} S_n(\cL^{-1}\bx_i, \wt\bp_m) f_2(\wt\bp_m),\quad \forall \bx_i \in X_1.
\end{equation}
\end{enumerate}
The evaluation of each $S_n(\cL^{-1}\bx_i, \wt\bp_m)$ or $S_n(\cR^{-1} \by, \wt\bq_{m'})$ has computational complexity $\cO(n^d)$. If the evaluation of $\cK(\cL\wt\bp_m, \cR\wt\bq_{m'})$ has only complexity $\cO(1)$, which happens when the coefficients are constants for instance, then the total complexity of evaluating~\eqref{EQ: SUMMATION} is $\cO(|X_2| n^{2d}) + \cO(n^{2d}) + \cO(|X_1| n^{2d}) \le \cO(N n^{2d})$, where $|X_i|$ is the number of collocation points in $X_i$ and $N$ is total number of collocation points in $\Omega$. Unfortunately, in general, the evaluation of $\cK(\cL\wt\bp_m, \cR\wt\bq_{m'})$ is not an order $\cO(1)$ operation for the variable coefficients case. We will discuss the cost for this situation in the next section.

\begin{remark}
    If all the cells $\{T_{j,h} \}_{j=1}^N$ are identical, then following Remark~\ref{REM: NEW KERNEL}, we can use $\wt\cK$ to replace $\cK$ in above Chebyshev interpolation formulation~\eqref{EQ:Kernel Approx} for higher accuracy. However, the evaluation of $\wt\cK$ costs much more than that of $\cK$ as we have discussed.
\end{remark}

%\begin{remark}
%The approximation error of~\eqref{EQ:Kernel Approx} is measured in $L^{\infty}$ sense, the convergence rate with respect to the order $n$ could be downgraded if the kernel function $\cK(\bx, \by)$'s regularity is limited. On the other hand, the approximation error can be also measured in $L^2$ sense, which is closely related to the singular-value-decomposition (SVD), which might relax the regularity assumption here.
%\end{remark}

%%%%%%%%%%%%%%%%%%%%%%%%%%%%%%%%%%%%%%%%%%%%%%%%%%%%%%%%%%%%%%%%%%
%%%%%%%%%%%%%%%%%%%%%%%%%%%%%%%%%%%%%%%%%%%%%%%%%%%%%%%%%%%%%%%%%%
\section{Implementation issues}
\label{SEC:Impl}
%%%%%%%%%%%%%%%%%%%%%%%%%%%%%%%%%%%%%%%%%%%%%%%%%%%%%%%%%%%%%%%%%%
%%%%%%%%%%%%%%%%%%%%%%%%%%%%%%%%%%%%%%%%%%%%%%%%%%%%%%%%%%%%%%%%%%

We now discuss briefly about some important issues on the implementation of the approach (M-iii) in Section~\ref{SEC:LINEAR SYS} with FMM algorithm described in the previous section.

%%%%%%%%%%%%%%%%%%%%%%%%%%%%%%%%%%%%%%%%%%%%%%%%%%%%%%%%%%%%%%%%%%
\subsection{Validity of low rank approximation} 
%%%%%%%%%%%%%%%%%%%%%%%%%%%%%%%%%%%%%%%%%%%%%%%%%%%%%%%%%%%%%%%%%%

When the coefficients $\mu$ and $\mu_s$ are sufficiently smooth as we have assumed, the boundedness of the exponential factor $E(\bx, \by)$ implies that our kernel $\cK(\bx, \by)$ in~\eqref{EQ:Kernel} admits at least the same low-rank approximation as the kernel $|\bx-\by|^{1-d}$ for far-field interaction. Such kernel has been well-studied in the fast multipole method community~\cite{BeGr-WMMEP97,ChGrRo-JCP99,ChJiMi-Book01,GrRo-AN97,MaRo-SIAM07,YiBiZo-JCP03}. This justifies the Chebyshev interpolation in~\eqref{EQ:Kernel Approx}.

%%%%%%%%%%%%%%%%%%%%%%%%%%%%%%%%%%%%%%%%%%%%%%%%%%%%%%%%%%%%%%%%%%
\subsection{Computational complexity of FMM}
\label{SEC: LINE INT}
%%%%%%%%%%%%%%%%%%%%%%%%%%%%%%%%%%%%%%%%%%%%%%%%%%%%%%%%%%%%%%%%%%

For the computational complexity, the most expensive step is the evaluation of the far-field interaction in~\eqref{EQ: M2L} where we have to evaluate the integral kernel $\cK$ for different Chebyshev point pairs. Each evaluation requires the computation of a line integral of the total absorption coefficient $\mu$ along the line that connects the Chebyshev points in far-field. When the total absorption coefficient $\mu$ is constant, this evaluation is trivial. However, when $\mu$ is not constant, this evaluation is quite expensive. In the following, we discuss the computational complexities for two practically important situations.

\begin{enumerate}
    \item When $\mu$ is sufficiently smooth and explicitly known for all $\bx\in\Omega$, then the line integral of $\mu$ could be evaluated analytically. 
    In this case, the computational complexity of evaluating $\cK(\bx, \by)$ will be $\cO(1)$, and the total complexity for the FMM is $\cO(N)$.
    
    \item When $\mu$ is only defined on collocation points, the line integral of $\mu$ cannot be evaluated directly. Assume that the line connecting $\bx$ and $\by$ passes through cells $T_{j_1, h},\dots, T_{j_m, h}$. Let $\bv = (\bx - \by)/ |\bx - \by|$ the unit direction along the line. We then have
    \begin{equation}\label{EQ: LINE INT}
    \int_{0}^{|\bx - \by|} \mu(\bx - s\bv) ds = \sum_{k=1}^m \int_{l_{k}}^{l_{k+1}} \mu(\bx - s\bv) ds,
    \end{equation}
    where the segment between $\bx - l_k \bv$ and $\bx-l_{k+1}\bv$ lies inside the cell $T_{j_k, h}$. Locally on each $T_{j_k, h}$, we check that the line integral
    \begin{multline}\label{EQ: INT ERROR}
    \int_{l_k}^{l_{k+1}} \mu(\bx - s\bv) ds = \int_0^{l_{k+1} - l_k} \mu(\bx - l_k\bv  - s\bv) ds =\\
    (l_{k+1} - l_k) \big( \mu(\bx_{j_k}) + \nabla \mu(\bx_{j_k})\cdot (\bx - l_k \bv -\bx_{j_k})\big) + \frac{1}{2}(l_{k+1} - l_k)^2 \bv\cdot \nabla\mu(\bx_{j_k}) + \cO(h^3).
    \end{multline}
    Using the fact that $|\bx - l_k \bv - \bx_{j_k}|\le h$ and $l_{k+1} - l_k \le h$, we conclude that we can simply approximate $\nabla \mu(\bx_{x_{j_k}})$ with first-order finite difference from neighboring collocation points, which leads to a truncation error of $\cO(h^3)$ according to ~\eqref{EQ: INT ERROR} (which means that the approximation error of the line integral~\eqref{EQ: LINE INT} is $\cO(h^2)$). This approximation is accurate enough to \emph{not} alter the order of error in~\eqref{EQ: ERROR}. However, the above approximation of line integral only uses local information and could lose smoothness across the cell boundaries where $\mu$'s value jumps. Therefore such formulation may limit the accuracy of FMM approximation when $\mu$ is not smooth.
\end{enumerate}

The overall complexity of the FMM algorithm follows directly from the discussion in~\cite{FoDa-JCP09}. Suppose the k-d tree for the algorithm has depth $L$, and at the $k$-th level there are $2^{kd}$ nodes $\cN_1, \dots, \cN_{2^{kd}}$. For each $\cN_i$, $1\le i\le 2^{kd}$, one has at most $(6^d - 3^d)$ nodes in its far-field. The corresponding evaluation cost is therefore at most $\cO(2^{-k}/h)$. Therefore, ignoring the constants in the $\cO$-notations, the total complexity will be dominated by:
\begin{multline}
    \sum_{k=1}^{L} 2^{kd} (6^d - 3^d) n^{2d} 2^{-k} h^{-1}  \le n^{2d} (6^d - 3^d) h^{-1} \frac{2^{(d-1)(L + 1)} - 1}{2^{d-1} - 1} \\  \le   n^{2d} (6^d - 3^d) h^{-1}\frac{2^{d-1} 2^{dL - L} }{2^{d-1} - 1} =\cO(2^{dL} 2^{-L}h^{-1})= \cO(N)
\end{multline}
where we have used the relation $N\simeq  n^d 2^{dL}$ and $h = \cO(N^{-1/d})$. This implies that the total complexity of the FMM is still at $\cO(N)$.

For the storage requirement, the usual $\cO(N)$ storage for the k-d tree structure is necessary. For a given order $n$ of the Chebyshev interpolation, we have $\cO(N)$ pairs of far field interactions for which we need to evaluate the kernel function. In our implementation, we cache all these kernel evaluations during the evaluation of $\phi$, which will be reused without any extra calculations during the \GMRES iterations; see, for instance, the numerical results in Table~\ref{TAB:Con Coeff} and Table~\ref{TAB:Var Coeff} of Section~\ref{SEC:Num}.

%In our im In the M2L step of the FMM algorithm for the evaluation of the sum~\eqref{EQ:Intgrl Discr}, there are at most $M = (7^d - 3^d)$ different transfer vectors $\bv\in \{\bv_1, \dots, \bv_M\}$ at top level. For each level $0\le k\le L$, the length of $t$ scales with the factor $2^{-k}$. Therefore, we can evaluate the integral in ~\eqref{EQ:FFT Mu} with for a cost on the order of $\cO(MQ)$. This means that the M2L step will have $\cO(Q n N)$ complexity.

%The second technique for accelerating the calculation is to use SVD compression. This works since the operator is compact and therefore its singular values decay rapidly. Under general setting, the kernel is not homogeneous nor transitional invariant, therefore pre-computation of SVD will cost $\mathcal{O}(N)$,  however it would still be helpful for GMRES iterations if singular values decay very rapidly.

%%%%%%%%%%%%%%%%%%%%%%%%%%%%%%%%%%%%%%%%%%%%%%%%%%%%%%%%%%%%%%%%%%
\subsection{Accuracy of the algorithm} 
%%%%%%%%%%%%%%%%%%%%%%%%%%%%%%%%%%%%%%%%%%%%%%%%%%%%%%%%%%%%%%%%%%

The accuracy of solution to the ERT~\eqref{EQ:ERT} with our numerical procedure depends mainly on two factors: the resolution of the spatial discretization $h$, and the accuracy of the fast multipole approximation of the summation~\eqref{EQ:Kernel Approx}. 
The latter relies on the order of the Chebyshev polynomial being used. Increasing the order of the polynomial should increase the accuracy of the approximation in general. However, that will also increase the computational cost of the algorithm, due to the increased cost in evaluating $S_n$. When $\mu$ increases significantly, the kernel function decays faster, hence it will require more points to resolve the kernel function. If the discretization is fixed, as $\mu$ increases, the accuracy can be reduced.

The convergence of the algorithm depends on the formulation of the matrix $K$. In Section~\ref{SEC:LINEAR SYS}, we have introduced a few methods based on the piecewise constant collocation method. One can construct higher ordered methods by a more careful treatment of the boundary layer effect. Due to the fact that an analytic form of solution to the integral equation~\eqref{EQ:ERT Intgrl U} can not be found on compact domains, we are only able to perform self-convergence test in Section~\ref{SEC:Num}.

%Let $U^{dir}$ be the numerical solution with a direct evaluation of the summation in~\eqref{EQ:Kernel Approx} and $U^{FMM}$ the FMM-accelerated numerical solution. When $U^{dir}$ and $U^{FMM}$ solutions are computed on the same mesh, the finer mesh will produce larger error $U^{dir}-U^{FMM}$ when the same number of interpolation points is used. This is because finer mesh provides structures that are harder to capture with the same interpolation polynomial. Moreover, the accuracy of approximating $U^{dir}$ by $U^{FMM}$ depends on the total absorption coefficient $\mu$ since the larger $\mu$ is, the faster the exponential decay is in the integral kernel. Therefore, for the same order of interpolation, the larger $\mu$ is, the worse the approximation is. We observe these phenomenon in our numerical experiments; see for instance, the simulations in Section~\ref{SEC:Num}.

%%%%%%%%%%%%%%%%%%%%%%%%%%%%%%%%%%%%%%%%%%%%%%%%%%%%%%%%%%%%%%%%%%
\subsection{Preconditioning techniques} 
\label{SEC: PRED}
%%%%%%%%%%%%%%%%%%%%%%%%%%%%%%%%%%%%%%%%%%%%%%%%%%%%%%%%%%%%%%%%%%

We developed two strategies for the preconditioning of the discrete linear system~\eqref{EQ: LINEAR SYSTEM} in special cases.

\paragraph{FFT-based Preconditioner.} The first preconditioning method works when the total absorption coefficient $\mu$ varies little on $\Omega$. In this case, the solution to the integral equation is very close to the solution to the same equation with a constant absorption coefficient $\overline{\mu}$, the mean value of $\mu$ over $\Omega$. We can therefore use the integral operator with $\overline{\mu}$ as a preconditioner for the true integral operator, since the integral operator with constant $\overline{\mu}$ is much cheaper to build. Moreover, we can use fast Fourier transform (FFT) techniques to accelerate the computation in domains of regular shapes. To be precise, let us write the integral formulation of the corresponding transport equation with $\overline{\mu}$ as
\begin{equation}
U - W\ast (\overline{\mu}_s U) = W \ast f
\end{equation}
where $W(\bx) = \exp(-\overline{\mu}|\bx|)/|\bx|^{d-1}$ and $*$ denotes the convolution product on $\Omega$. Without loss of generality, let us assume that $\Omega\subset [0, 1]^d$ and extend the convolution kernel $W$ to a periodic function $\widetilde{W}$ on $[0, 2]^d$. We take the zero extensions of $\overline{\mu}_s U$ and $f$ as $\widetilde{\overline{\mu}_s U}$ and $\widetilde{f}$ from $\Omega$ to $[0, 2]^d$ respectively. The above equation is transformed to
\begin{equation}\label{EQ: FFT}
P \cF^{-1}\left(\cF\left(\overline\mu_s^{-1}\chi_{[0, 2]^d} -  \widetilde{W}\right)\cF (\widetilde{\overline\mu_s U}) \right) = W\ast f
\end{equation}
where $\cF$ denotes the Fourier transform, $P$ is the restriction operator from $[0, 2]^{d}$ back to $\Omega$, and $\chi_{[0, 2]^d}$ is the characteristic function of ${[0, 2]^d}$. Since $P: [0,2]^d\mapsto \Omega$ is not invertible, the above equation cannot be solved directly. Therefore we try to regularize the operator $P$. Here we make the additional assumption that the support of source function $f$ is far away from the boundary $\partial\Omega$. 
With these assumptions, we can simply regularize the inversion of $P$ by taking zero padding,
\begin{equation}
	P^{-1} g = g\chi_{\Omega}.
\end{equation}
%Physically speaking, since $\overline\mu_s$ is relatively large in practice, when the above assumptions are fulfilled, the source function will be mostly absorbed and only a small amount of radiation will leak to the outside. 
Therefore we can approximately solve~\eqref{EQ: FFT} by
\begin{equation}\label{EQ:SOL FFT}
U \simeq \frac{1}{\overline\mu_s} P\cF^{-1} \left( \cF \left( P^{-1}(W\ast f) \right) / \cF(\overline\mu_s^{-1}\chi_{[0, 2]^d} - \widetilde{W}) \right).
\end{equation}

\paragraph{Diffusion Preconditioner.} The second preconditioning method we implemented is based on the diffusion approximation (DA) of the radiative transport equation. Diffusion approximation is valid when the characteristic scale of the domain, say $\ell_{\Omega} := \text{diam}(\Omega)$, is very large compared to the meas free path of the transport problem, $\sim (d\mu_s)^{-1}$, that is when $d\mu_s \ell_{\Omega}$ is very large. In that case, it can be shown that the solution to the integral equation~\eqref{EQ:ERT Intgrl U} solves approximately the following diffusion equation~\cite{DaLi-Book93-6}:
\begin{equation}\label{EQ:DA}
\begin{array}{rcll}
-\nabla \cdot \dfrac{1}{d \mu(\bx)} \nabla U + \mu_a U &=& f(\bx), &\text{in}\ \ \Omega \\
U + \ell \dfrac{\partial U}{\partial \bn} &=& 0, &\text{on}\ \partial\Omega
\end{array}
\end{equation} 
where $\mu_a=\mu-\mu_s$ is the physical absorption coefficient and $\ell\simeq \frac{2}{d\mu}$ is called the extrapolation length~\cite{DaLi-Book93-6}. 

In diffusive regime, the solution of~\eqref{EQ:DA}, which is computationally inexpensive to obtain, is a good approximation to the solution of ~\eqref{EQ:ERT Intgrl U}. This is the main motivation for the diffusion operator as a preconditioner for the transport problem, leading to the popular diffusion synthetic acceleration (DSA) technique~\cite{Alcouffe-NSE77,Larsen-TTSP84}. Our implementation of the diffusion preconditioner is, however, more natural than the classical DSA implementation since we work directly on the variable $U$, the integral of the transport solution over the angular variable, not the transport solution itself. In our implementation of the diffusion preconditioner, we discretize~\eqref{EQ:DA} with a classical second order finite difference method, and apply the diffusion solver with a standard V-cycle multigrid method. The transfer between the diffusion solution and the solution to the integral equation~\eqref{EQ:ERT Intgrl U} are done with a standard intepolation algorithm.

%%%%%%%%%%%%%%%%%%%%%%%%%%%%%%%%%%%%%%%%%%%%%%%%%%%%%%%%%%%%%%%%%%
%%%%%%%%%%%%%%%%%%%%%%%%%%%%%%%%%%%%%%%%%%%%%%%%%%%%%%%%%%%%%%%%%%
\section{Numerical experiments}
\label{SEC:Num}
%%%%%%%%%%%%%%%%%%%%%%%%%%%%%%%%%%%%%%%%%%%%%%%%%%%%%%%%%%%%%%%%%%
%%%%%%%%%%%%%%%%%%%%%%%%%%%%%%%%%%%%%%%%%%%%%%%%%%%%%%%%%%%%%%%%%%

We now present some numerical simulations to demonstrate the performance of the algorithm we developed. We perform simulations in both homogeneous and inhomogeneous media. In the homogeneous case, both $\mu_a$ and $\mu_s$ are constants (and so is $\mu=\mu_a+\mu_s$). We can therefore evaluate $E(\bx, \by)= \exp(-\mu |\bx - \by|)$ analytically. In the inhomogeneous case, both $\mu_a$ and $\mu_s$ are at least $C^2$ smooth but are defined on the collocation points only. We therefore have to compute the line integral in $E(\bx, \by)$ by the summation of line integrals on each cell; see the discussions in Section~\ref{SEC: LINE INT} on this issue.
\begin{figure}[!htb]
	\centering
	\includegraphics[width=0.4\textwidth,height=0.33\textwidth]{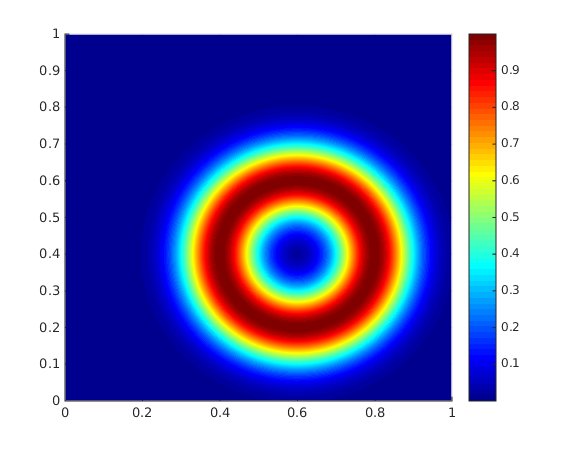} \hskip 0.5cm
	\includegraphics[width=0.4\textwidth,height=0.33\textwidth]{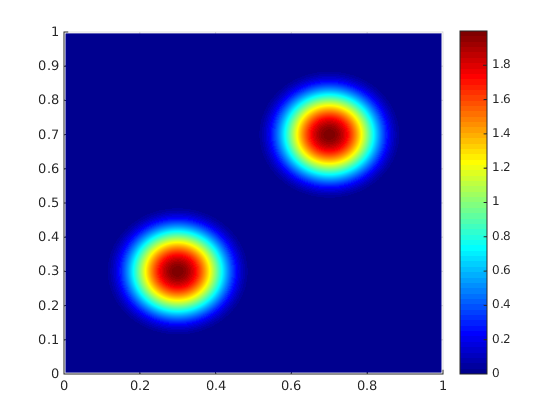}
	\caption{The two smooth source functions used in the numerical experiments.}
	\label{FIG:Sources}
\end{figure}

All the simulations are performed on the nondimensionalized transport equation so that the physics are on determined by the relative size of the coefficients to the size of the domain of the problem which we take as the unit square $\Omega=(0, 1)^2$. We discretize the domain with uniform grids as described in Section~\ref{SEC:DISC}. The linear systems involved are also solved with the \GMRES algorithm implemented in MATLAB. We set the \GMRES stopping tolerance at $10^{-12}$. We vary the scattering coefficient $\mu_s$ to test the performances of the algorithm in different regimes.  The larger the scattering coefficient $\mu_s$ is, the more diffusive the solution behaves. As we will see in the following experiments, the performance of our algorithm does not change dramatically from the low scattering transport regime to the moderately scattering transport regime.

To simplify the presentation, we list in Table~\ref{TAB: NOTATION} various parameters to be used to measure the performance of the algorithms in solving the integral equation~\eqref{EQ:ERT Intgrl FMM}. Note that in our computations, we have cached all the line integrals needed when setting up the algorithm. Therefore, the cost $T^{\text{\tiny{fmm}}}_{\text{\tiny{g}}}$ (resp. $T^{\text{\tiny{dir}}}_{\text{\tiny{g}}}$) does not include the cost $T^{\text{\tiny{fmm}}}$ (resp. $T^{\text{\tiny{dir}}}$). All the computations shown below are \emph{single-threaded} on a Linux laptop with i7-6560U CPU @ 2.20GHz and 16GB RAM. 
\begin{table}[!htb]
\centering
\scriptsize
\caption{List of parameters used to measure the performance of the algorithms to solve~\eqref{EQ:ERT Intgrl FMM}}
\vspace{0.2cm}
\label{TAB: NOTATION}
\begin{tabular}{|l|l|}
\hline 
$N$ & total number of collocation points \\
\hline 
$n$ & order of Chebyshev interpolation; see~\eqref{EQ: CHEB}\\
\hline
$K$ & kernel matrix $K=[K_{jk}\omega_k]$; see~\eqref{EQ: LINEAR SYSTEM}\\
\hline
$\Td$ & time cost (seconds) of direct matvec of $K$ \\
\hline
$\Tf$ & time cost (seconds) of FMM accelerated matvec of $K$ \\
\hline
$\Tdg$ & time cost (seconds) of each iteration in \GMRES with direct matvec\\
\hline
$\Tfg$ & time cost (seconds) of each iteration in \GMRES with FMM acceleration without preconditioning\\
\hline
$\Tfpg$ &  time cost (seconds) of each iteration in \GMRES FMM acceleration with preconditioner\\
\hline
$\Idg$ & total iteration number of \GMRES with direct matvec\\
\hline
$\Ifg$ & total iteration number of \GMRES with FMM acceleration without preconditioning\\
\hline
$\Ifpg$ & total iteration number of \GMRES with FMM acceleration with preconditioner\\
\hline
$E_{\ell^2}$ & relative $\ell^2$ error between the solutions from direct matvec and FMM acceleration. \\
\hline
%$\text{E.O.C}$ & estimated order of convergence\\
%\hline
\end{tabular}
\end{table}

%%%%%%%%%%%%%%%%%%%%%%%%%%%%%%%%%%%%%%%%%%%%%%%%%%%%%%%%%%%%%%%%%%
\subsection{Validation of low rank property}
\label{EXP:I}
%%%%%%%%%%%%%%%%%%%%%%%%%%%%%%%%%%%%%%%%%%%%%%%%%%%%%%%%%%%%%%%%%%

We first validate the low rank property for the matrix $K$ in~\eqref{EQ: LINEAR SYSTEM} by following the approach (M-iii) in the discretization; see Section~\ref{SEC:LINEAR SYS}. Let $U^{\rm fmm}$ and $U^{\rm dir}$ be the solutions to~\eqref{EQ: LINEAR SYSTEM} with and without using FMM to accelerate respectively. We use the relative $\ell^2$ error $$E_{\ell^2}=\frac{ \|U^{\rm dir}-U^{\rm fmm}\|_{l^2}}{\|U^{\rm dir}\|_{\ell^2}}$$ to measure their difference.

\paragraph{Experiment I.} In the first numerical experiment, we perform simulations with a constant scattering coefficient $\mu_s\equiv 2.0$ and a constant total absorption coefficient $\mu\equiv 2.2$ (which means the physical absorption is $\mu_a \equiv 0.2$). The source function $f(\bx)$ is the ring source illustrated in the left plot of Figure~\ref{FIG:Sources}. In Table~\ref{TAB:Con Coeff} we show comparisons in three groups with increasing order of Chebyshev interpolation: $n=4$, $n=6$ and $n=9$. We first note that, with reasonable relative approximation error $E_{\ell^2}$ (e.g. on the order of $10^{-4}$ with $n=4$), the growth of running time with respect to $N$ is almost linear for FMM accelerated \GMRES (see $\Tf$, $\Tfg$) and quadratic for the regular \GMRES (see $\Td$, $\Tdg$). When $N$ is relatively small, due to the fast DGEMM implementation in BLAS called from MATLAB, the direct matrix-vector production is slightly faster than the FMM acceleration, while for larger $N$, the FMM accelerated \GMRES outperforms the regular \GMRES. This trend is kept when we increase the accuracy of the FMM approximation by increasing, $n$, the order of Chebyshev interpolation. When the spatial discretization gets too fine (for instance when $N\ge 65536$), it takes the regular \GMRES algorithm too much memory and time to finish the calculations. We have to stop the algorithm before it converges. However, the FMM accelerated \GMRES can still solve the linear system~\eqref{EQ: LINEAR SYSTEM} in relatively short time. In Table~\ref{TAB:Con Coeff}, when $N=1024$ and $n=9$, the relative $\ell^2$ error is comparable to machine precision, since the FMM here only contains near-field interactions.
\begin{table}[!htb]
\centering
\small
\caption{The computational costs and relative errors between the solutions with and without FMM acceleration for a homogeneous media with $\mu_s=2.0$ and $\mu_a=0.2$ under various total collocation points $N$ and Chebyshev interpolation orders $n$.}
\label{TAB:Con Coeff}
    \begin{tabular}{l*{10}{c}l}
        \hline
        $N$  & $n$  & $\Tf$ & $\Tfg$ & $\Ifg$ & $\Td$ & $\Tdg$ & $\Idg$ & $E_{\ell^2}$ \rule{0pt}{2.6ex}\rule[-1.2ex]{0pt}{0pt} \\
        \hline
        \rule{0pt}{4ex}1,024 &  4  & $\ep{2.11}{-2}$ & $\ep{5.83}{-3}$ & 10 & $\ep{6.76}{-2}$ & $\ep{4.20}{-4}$& $10$ & $\ep{8.53}{-5}$ \\
        4,096 & 4  & $\ep{1.06}{-1}$ & $\ep{1.02}{-2}$ & 10 & $\ep{1.08}{+0}$  & $\ep{7.94}{-3}$& 10& $\ep{1.12}{-4}$ \\
        16,384 & 4  & $\ep{4.49}{-1}$ & $\ep{5.48}{-2}$ & 10 & $\ep{1.90}{+1}$ & $\ep{1.21}{-1}$  & 10 &$\ep{1.22}{-4}$\\
        65,536 & 4  & $\ep{1.89}{+0}$ & $\ep{2.69}{-1}$ & 10 & -- & --&--&-- \\
        262,144 & 4  & $\ep{8.82}{+0}$  & $\ep{2.46}{+0}$ &10 & -- & --&--&-- \\
        \\
        1,024 & 6 &$\ep{2.21}{-2}$& $\ep{2.22}{-3}$ &10 & $\ep{6.76}{-2}$ & $\ep{4.20}{-4}$&  10 &$\ep{1.02}{-6}$\\
        4,096 & 6 &$\ep{1.37}{-1}$& $\ep{1.57}{-2}$& 10 & $\ep{1.08}{+0}$ & $\ep{7.94}{-3}$& 10 & $\ep{1.28}{-6}$\\
        16,384 & 6 & $\ep{6.74}{-1}$& $\ep{6.42}{-2}$&10 &$\ep{1.90}{+1}$& $\ep{1.21}{-1}$ & 10 & $\ep{1.31}{-6}$\\
        65,536 & 6  & $\ep{2.96}{+0}$ & $\ep{3.29}{-1}$ & 10 &--& --&--&--\\
        262,144 & 6  & $\ep{1.39}{+1}$  & $\ep{2.71}{+0}$ &10 & -- & --&--&-- \\
        \\
        1,024& 9 & $\ep{6.52}{-2}$& $\ep{3.06}{-3}$ &10 & $\ep{6.76}{-2}$ & $\ep{4.20}{-4}$&  10& $\ep{5.70}{-16}$\\
        4,096& 9 & $\ep{2.53}{-1}$& $\ep{1.69}{-2}$ &10& $\ep{1.08}{+0}$ & $\ep{7.94}{-3}$& 10&  $\ep{3.16}{-9}$ \\
        16,384&9 & $\ep{1.36}{+0}$& $\ep{9.25}{-2}$ &10& $\ep{1.90}{+1}$ & $\ep{1.21}{-1}$ &10& $\ep{2.45}{-9}$\\
        65,536 & 9  & $\ep{6.42}{+0}$ & $\ep{4.88}{-1}$ & 10 & -- & --&--&-- \\
        262,144 & 9  & $\ep{3.06}{+1}$  & $\ep{3.40}{+0}$ &10 & -- & --&--&--\\
    \end{tabular}
\end{table}
\begin{table}[!htb]
    \centering
    \small
    \caption{The computational costs and relative errors between the solutions with and without FMM acceleration for for the heterogeneous medium in~\eqref{EQ:Coeff} under various total collocation points $N$ and Chebyshev interpolation orders $n$.}
    \label{TAB:Var Coeff}
    \begin{tabular}{l*{10}{c}l}
        \hline
        $N$  & $n$  & $\Tf$ & $\Tfg$ & $\Ifg$ & $\Td$ & $\Tdg$ & $\Idg$ & $E_{\ell^2}$ \rule{0pt}{2.6ex}\rule[-1.2ex]{0pt}{0pt} \\
        \hline
        \rule{0pt}{4ex}1,024 &  4  & $\ep{2.35}{-1}$ & $\ep{4.18}{-3}$ & 15 & $\ep{1.34}{+1}$ & $\ep{5.98}{-4}$& $15$ & $\ep{2.00}{-4}$ \\
        4,096 & 4  & $\ep{1.20}{+0}$ & $\ep{7.50}{-3}$ & 15 & $\ep{2.20}{+2}$  & $\ep{7.42}{-3}$& 15& $\ep{3.07}{-4}$ \\
        16,384 & 4  & $\ep{5.95}{+0}$ & $\ep{3.61}{-2}$ & 15 & $\ep{3.73}{+3}$ & $\ep{1.15}{-1}$  & 15 &$\ep{3.54}{-4}$\\
        65,536 & 4  & $\ep{2.67}{+1}$ & $\ep{2.32}{-1}$ & 15 & -- & --&--&-- \\
        262,144 & 4  & $\ep{1.14}{+2}$  & $\ep{2.22}{+0}$ &15 & -- & --&--&-- \\
        \\
        1,024 & 6 &$\ep{3.66}{-1}$& $\ep{4.24}{-3}$ &15 & $\ep{1.34}{+1}$ & $\ep{5.98}{-4}$& $15$ &$\ep{1.73}{-5}$\\
        4,096 & 6 &$\ep{2.62}{+0}$& $\ep{9.56}{-3}$& 15 & $\ep{2.20}{+2}$ & $\ep{7.42}{-3}$& 15 & $\ep{1.37}{-5}$\\
        16,384 & 6 & $\ep{1.45}{+1}$& $\ep{4.66}{-2}$&15 &$\ep{3.73}{+3}$& $\ep{1.15}{-1}$ & 15 & $\ep{7.05}{-6}$\\
        65,536 & 6  & $\ep{7.01}{+1}$ & $\ep{2.69}{-1}$ & 15 &--& --&--&--\\
        262,144 & 6  & $\ep{3.08}{+2}$  & $\ep{2.39}{+0}$ &15 & -- & --&--&-- \\
        \\
        1,024& 9 & $\ep{7.71}{-1}$& $\ep{2.57}{-2}$ &15 & $\ep{1.34}{+1}$ & $\ep{5.98}{-4}$& $15$ & $\ep{6.99}{-16}$\\
        4,096& 9 & $\ep{7.15}{+0}$& $\ep{1.49}{-2}$ &15& $\ep{2.20}{+2}$ & $\ep{7.42}{-3}$& 15&  $\ep{4.94}{-6}$ \\
        16,384&9 & $\ep{4.42}{+1}$& $\ep{7.52}{-2}$ &15& $\ep{3.73}{+3}$ & $\ep{1.15}{-1}$ &15& $\ep{3.03}{-6}$\\
        65,536 & 9  & $\ep{2.24}{+2}$ & $\ep{3.94}{-1}$ & 15 & -- & --&--&-- \\
        262,144 & 9  & $\ep{1.10}{+3}$  & $\ep{3.29}{+0}$ &15 & -- & --&--&--\\
    \end{tabular}
\end{table}

\paragraph{Experiment II.} In the second numerical experiment, we repeat the simulations in Experiment I for a heterogeneous medium. The coefficients are given as
\begin{equation}\label{EQ:Coeff}
	\mu_a(\bx) = 0.2, \qquad 
	\mu_s(\bx) = 3.0 + 2.0  \exp\left(-\frac{(x - 0.5)^2 + (y - 0.5)^2}{4}\right).
\end{equation}
Our algorithm only uses the values of the coefficients on the collocation points. We again use the the ring source illustrated in the left plot of Figure~\ref{FIG:Sources}. In Table~\ref{TAB:Var Coeff} we show comparisons in three groups with increasing order of Chebyshev interpolation. The first noticeable difference between Table~\ref{TAB:Con Coeff} and Table~\ref{TAB:Var Coeff} is that the computational cost to evaluate the matrix-vector multiplication is now considerably higher. This is mainly due to the fact that for the variable coefficients case, we need to evaluate the line integrals by summation of quadratures in each cell~\eqref{EQ: INT ERROR}, while in the constant coefficients case the kernels are given analytically for any pair $(\bx, \by)$. In our implementation, we cached all the line integrals so that they can be used repeatedly during \GMRES iterations. This is the reason why the time cost in each iteration $\Tfg$ for variable coefficient cases in Table~\ref{TAB:Var Coeff} is very similar to the corresponding constant coefficient cases in Table~\ref{TAB:Con Coeff}. The overall computational costs again scale almost linearly with respect to the total collocation points $N$. Another noticeable difference between results in Table~\ref{TAB:Con Coeff} and Table~\ref{TAB:Var Coeff} is that, when we increase the order of Chebyshev interpolation $n$, the relative error $E_{\ell^2}$ of Table~\ref{TAB:Var Coeff} does not decay as fast as the relative error in Table~\ref{TAB:Con Coeff}. This is caused by the in accuracy of $E(\bx, \by)$ across the boundaries of the grid cells resulted from the local approximation of the line integrals in~\eqref{EQ: INT ERROR}. Similar to Experiment I,  when $N=1024$ and $n=9$, the relative $\ell^2$ error is comparable to machine precision for the same reason.

%%%%%%%%%%%%%%%%%%%%%%%%%%%%%%%%%%%%%%%%%%%%%%%%%%%%%%%%%%%%%%%%%%
\subsection{Accuracy of FMM solution}
\label{SUBSEC:Accuracy}
%%%%%%%%%%%%%%%%%%%%%%%%%%%%%%%%%%%%%%%%%%%%%%%%%%%%%%%%%%%%%%%%%%

We now present some self-convergence tests on the accuracy of the solutions to the integral equation~\eqref{EQ: LINEAR SYSTEM} by the FMM accelerated \GMRES algorithm.

\paragraph{Experiment III.} In this experiment, we perform simulations with constant scattering coefficients $\mu_s = 2.0,\, 5.0,\, 10.0$ and a fixed constant absorption coefficient $\mu_a\equiv 0.2$. The source function $f(\bx)$ is the left plot of Figure~\ref{FIG:Sources}. We choose the discretization as a uniform grid with cell size of $h = \frac{1}{24 (2k-1)}$ for $1\le k\le 12$ and take the solution of $k=12$ as the reference solution. The numerical errors are all evaluated on the common collocation points with $\ell^2$ norm on the coarsest level of $k=1$. We perform simulations based on approaches (M-ii) and (M-iii) in Section~\ref{SEC:LINEAR SYS} for the construction of the integral operator $K$. The order of Chebyshev interpolation is fixed as $n = 6$ for both simulations. 

\begin{figure}[!htb]
    \centering
    \includegraphics[height=200pt, width=220pt]{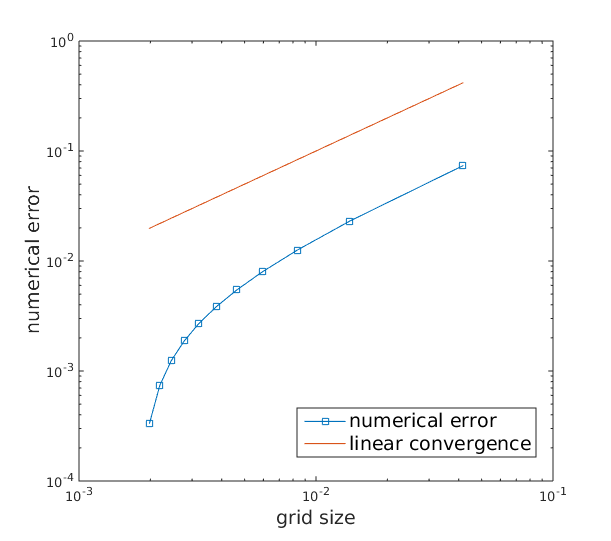}
    \includegraphics[height=200pt, width=220pt]{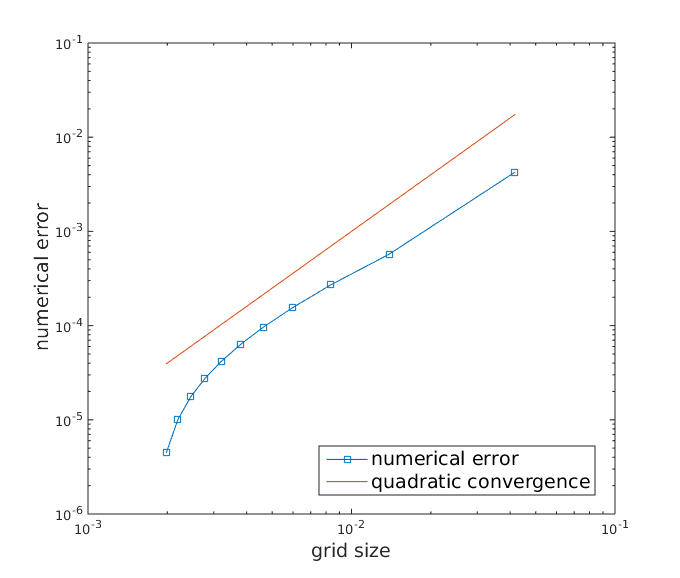}
    \caption{The numerical errors with respect to various grid sizes for $\mu_s = 2.0$ with the source function given on the left plot of Figure~\ref{FIG:Sources}. Shown are the $\ell^2$ error of the solutions compared with the reference solution calculated at $k=12$ using approaches (M-iii) (left) and (M-ii) (right) respectively.}
    \label{FIG: CONV1-1}
\end{figure}

\begin{figure}[!htb]
    \centering
    \includegraphics[height=200pt, width=220pt]{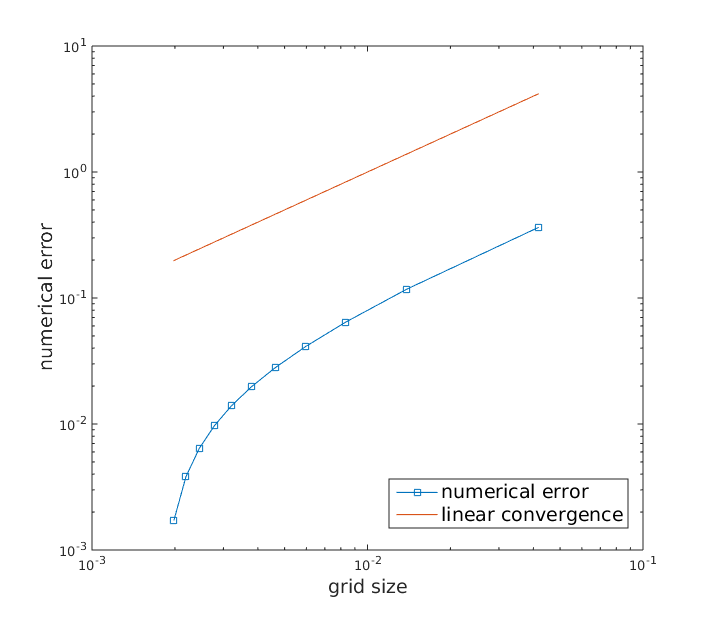}
    \includegraphics[height=200pt, width=220pt]{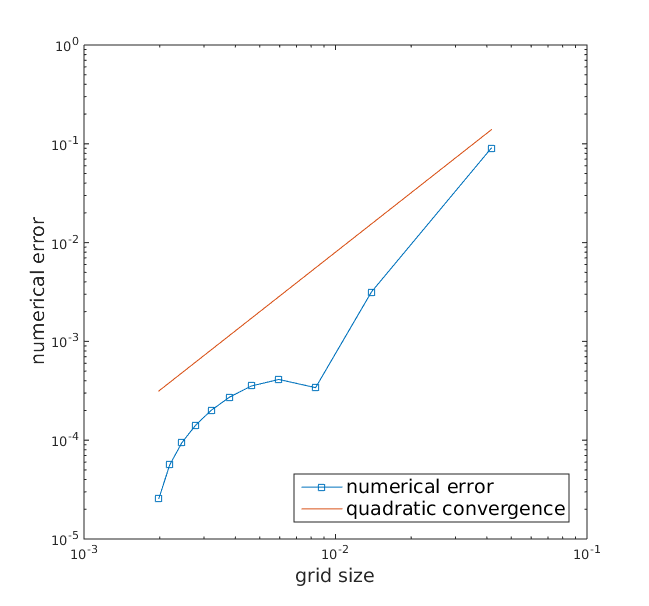}
    \caption{Same as Figure~\ref{FIG: CONV1-1} except that $\mu_s = 5.0$ here.}
    \label{FIG: CONV1-2}
\end{figure}

\begin{figure}[!htb]
    \centering
    \includegraphics[height=200pt, width=220pt]{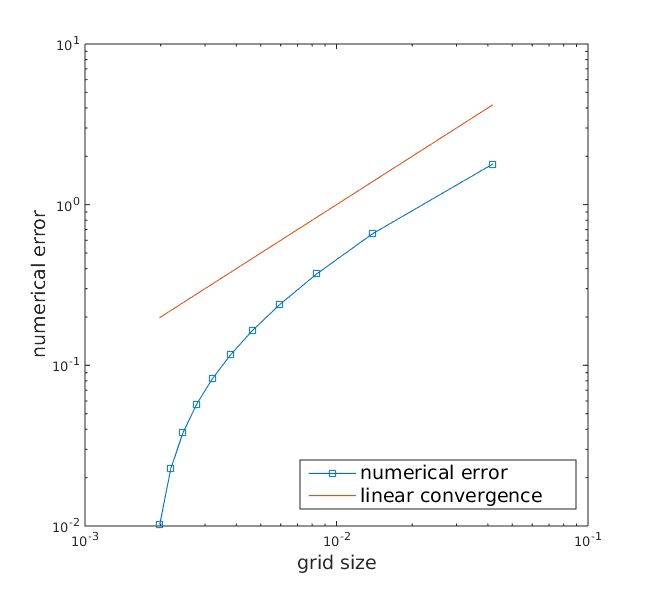}
    \includegraphics[height=200pt, width=220pt]{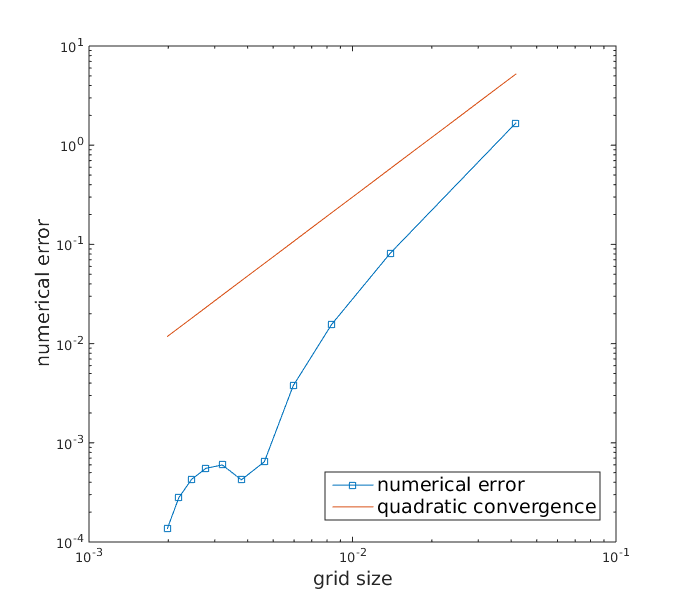}
    \caption{Same as Figure~\ref{FIG: CONV1-1} except that $\mu_s = 10.0$ here.}
    \label{FIG: CONV1-3}
\end{figure}

The numerical results for $\mu_s =2.0,\,5.0,\,10.0$ are shown  Figure~\ref{FIG: CONV1-1}, Figure~\ref{FIG: CONV1-2} and Figure~\ref{FIG: CONV1-3} respectively. The left and the right plots in each figure are results for approaches (M-iii) and (M-ii) respectively. We observe that, when the grid is relatively coarse, the convergence rates for approaches (M-iii) and (M-ii) are roughly linear and quadratic respectively. When the grids get very fine, better convergence behavior emerges. This is because in such cases, the evaluation points are far away from the boundary, which are therefore less affected by the boundary effect.

\begin{figure}[!htb]
    \centering
    \includegraphics[height=200pt, width=220pt]{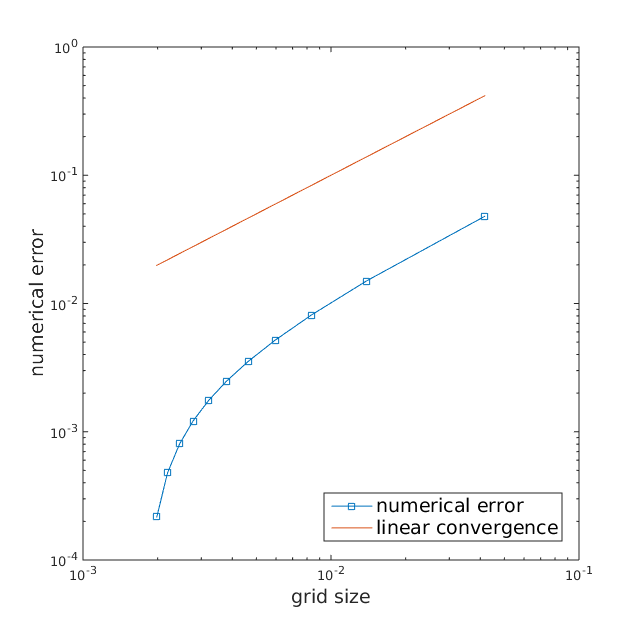}
    \includegraphics[height=200pt, width=220pt]{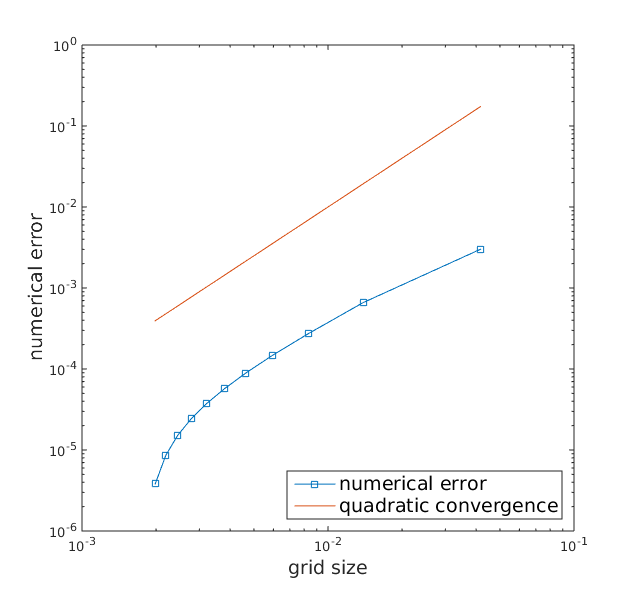}
    \caption{The numerical errors with respect to various grid sizes for $\mu_s = 2.0$ with the source function given on the right plot of Figure~\ref{FIG:Sources}. Shown are the $\ell^2$ error of the solutions compared with the reference solution calculated at $k=12$ using approaches (M-iii) (left) and (M-ii) (right) respectively.}
    \label{FIG: CONV2-1}
\end{figure}

\begin{figure}[!htb]
    \centering
    \includegraphics[height=200pt, width=220pt]{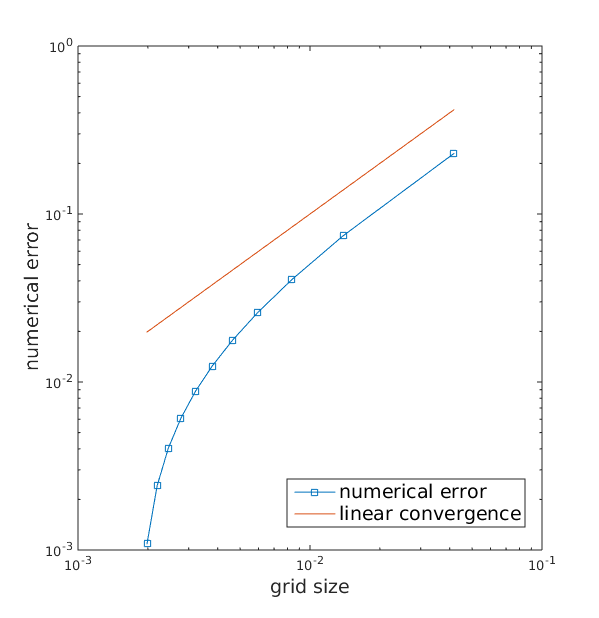}
    \includegraphics[height=200pt, width=220pt]{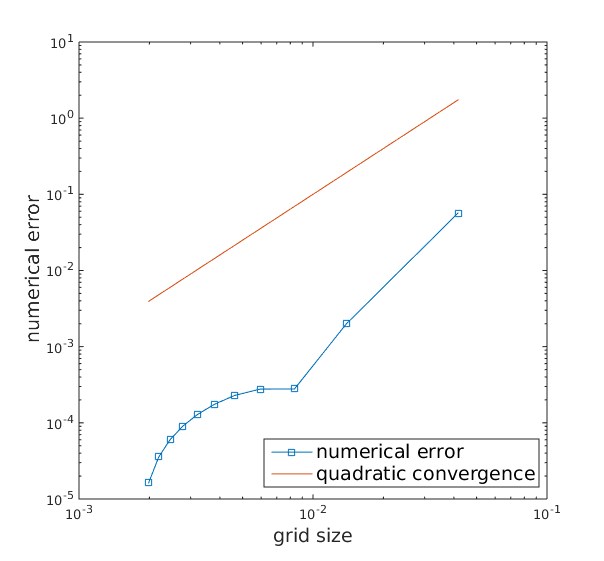}
    \caption{Same as Figure~\ref{FIG: CONV2-1} except that $\mu_s = 5.0$ here.}
    \label{FIG: CONV2-2}
\end{figure}

\begin{figure}[!htb]
    \centering
    \includegraphics[height=200pt, width=220pt]{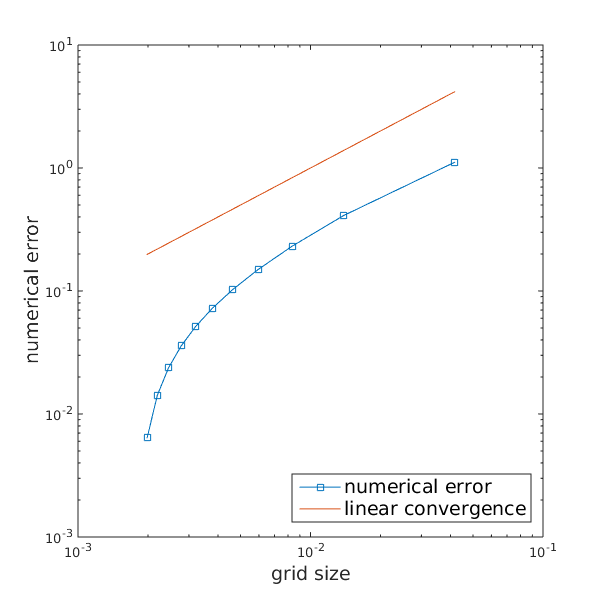}
    \includegraphics[height=200pt, width=220pt]{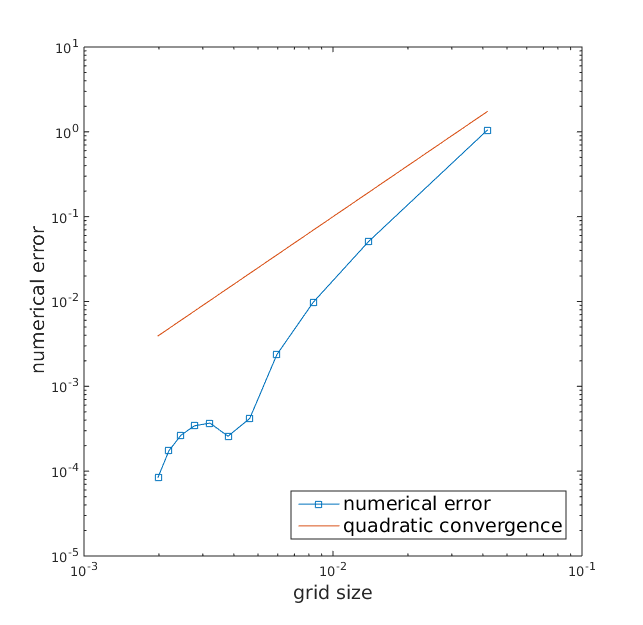}
    \caption{Same as Figure~\ref{FIG: CONV2-1} except that $\mu_s = 10.0$ here.}
    \label{FIG: CONV2-3}
\end{figure}

\paragraph{Experiment IV.} We repeat the numerical simulations in Experiment III for a different source function, given in the right plot of Figure~\ref{FIG:Sources}. The numerical results for $\mu_s =2.0,\,5.0,\,10.0$ are shown Figure~\ref{FIG: CONV2-1}, Figure~\ref{FIG: CONV2-2} and Figure~\ref{FIG: CONV2-3} respectively. What we observed in Experiment III can also be observed here: when the grids are relatively coarse, the convergence rates are almost linear for approach (M-iii) and quadrature for approach (M-ii), and when the grids get finer, better convergence behaviors emerge for both approaches. Moreover, in both Experiment III and Experiment IV, the numerical errors are larger when the scattering coefficient $\mu_s$ gets larger. This is mainly due to the fact that larger $\mu_s$ requires finer discretization for the kernel to capture the faster decay of the kernel induced by the fact $E(\bx, \by)$ which depends on $\mu_s$.

%%%%%%%%%%%%%%%%%%%%%%%%%%%%%%%%%%%%%%%%%%%%%%%%%%%%%%%%%%%%%%%%%%
\subsection{Preconditioning}
\label{SUBSEC:Preconditioning}
%%%%%%%%%%%%%%%%%%%%%%%%%%%%%%%%%%%%%%%%%%%%%%%%%%%%%%%%%%%%%%%%%%

We now study the performances of the two preconditioners we described in Section~\ref{SEC: PRED}. 

\begin{figure}[!htb]
    \centering
    \includegraphics[height=200pt, width=220pt]{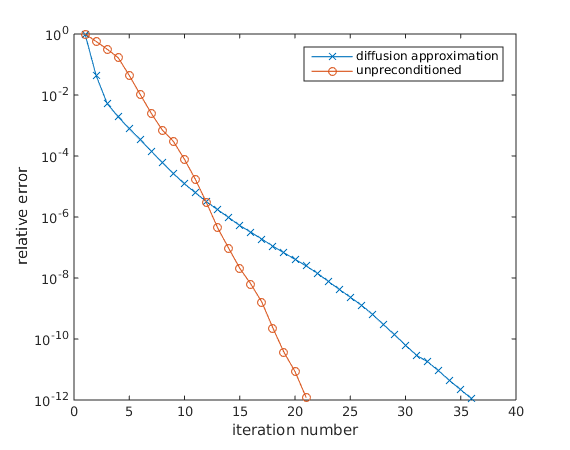}
    \includegraphics[height=200pt, width=220pt]{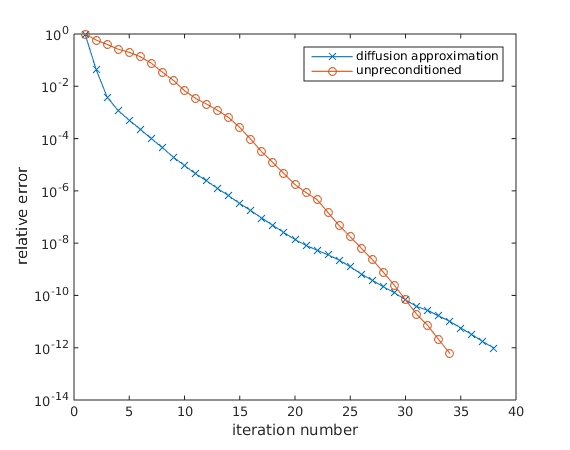}
    \includegraphics[height=200pt, width=220pt]{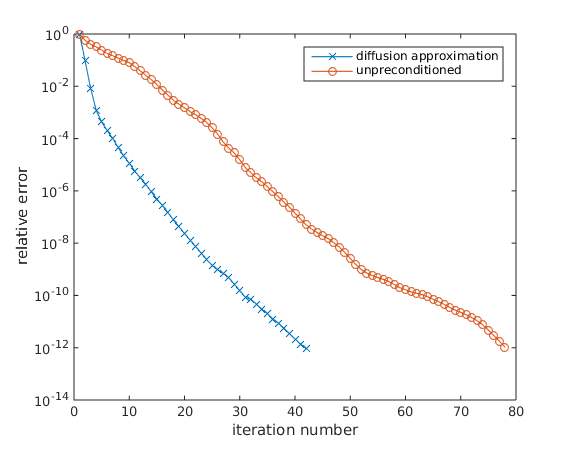}
    \includegraphics[height=200pt, width=220pt]{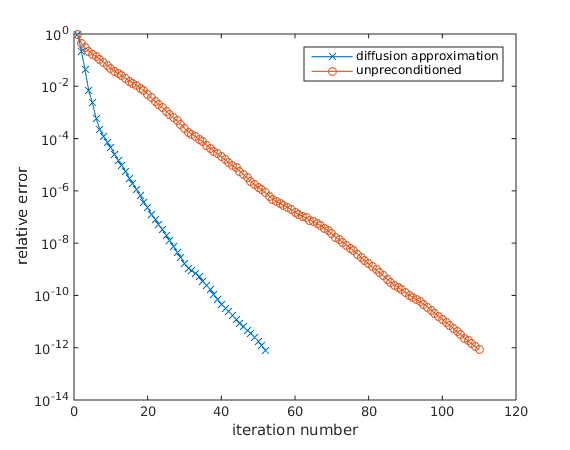}
    \caption{Convergence history for the diffusion preconditioned (blue $\times$) and unpreconditioned (orange $\circ$) \GMRES iterations with different scattering coefficients. Show from left and right are cases at $\mu_s=10.0$, $\mu_s = 20.0$, $\mu_s = 40.0$, and $\mu_s = 80.0$ respectively.}
    \label{FIG: EXP V}
\end{figure}

\paragraph{Experiment V.} We fix the constant absorption coefficient $\mu_a = 0.2$ and vary the scattering coefficient to values $\mu_s = 10.0,\, 20.0,\, 40.0,\, 80.0$. We take a grid with cell size $h = \frac{1}{512}$, which results in a total collocation number $N = 262144$ . The order of Chebyshev interpolation is fixed as $n=6$. We use the source function on the right plot of Figure~\ref{FIG:Sources}. The numerical results are shown in Figure~\ref{FIG: EXP V} and the Table~\ref{TAB:EXP VI}. The number of \GMRES iterations to reach the convergence tolerence we set, $\eps=10^{-12}$, for the preconditioned and unpreconditioned algorithms are respectively $36$ and $21$ ($\mu_s = 10.0$), $38$ and $34$ ($\mu_s=20.0$), $42$ and $78$ ($\mu_s=40.0$), and $52$ and $110$ ($\mu_s=80.0$). This shows that when $\mu_s$ is sufficiently large, say $\mu_s>20.0$, the diffusion preconditioner is effective in bringing down the total number of iterations it takes for the \GMRES algorithm to converge. However, our implementation of the diffusion preconditioner are not effective enough, so that the computational cost per iteration for the preconditioned version of the \GMRES iteration is significantly higher than the unpreconditioned iterations. Due to this, the overall cost of the preconditioned method (up to the convergence) is only lower than that of the unpreconditioned method when $\mu_s$ is very large (when $\mu_s>20.0$).   
\begin{table}[!htb]
     \centering
     \small
     \caption{The computational cost per iteration and the number of iterations to convergence for diffusion preconditioned and unpreconditioned \GMRES iterations under different scattering coefficients.}
     \vspace{0.2cm}
     \label{TAB:EXP V}   
         \begin{tabular}{l*{10}{c}l}
             \hline
             $\mu_s$  & $\Tfg$ & $\Ifg$ & $\Tfpg$ & $\Ifpg$ \rule{0pt}{2.6ex}\rule[-1.2ex]{0pt}{0pt} \\
             \hline
             \rule{0pt}{4ex}$10.0$ &  $\ep{2.50}{+0}$ & $21$ & $\ep{5.18}{+0}$ & $36$ \\
             $20.0$ & $\ep{3.13}{+0}$  & $34$ & $\ep{5.13}{+0}$ & $38$ \\
             $40.0$ & $\ep{2.95}{+0}$  & $78$ & $\ep{5.07}{+0}$ & $42$\\
             $80.0$ & $\ep{3.11}{+0}$  & $110$ & $\ep{5.07}{+0}$ & $52$  %\\
             %$\eqref{EQ:VAR MUS}$ & $\ep{3.06}{+0}$  & $113$ & $\ep{4.14}{+0}$ & $52$ 
            \end{tabular}
\end{table}
   
\begin{figure}[!htb]
    \centering
    \includegraphics[height=200pt, width=220pt]{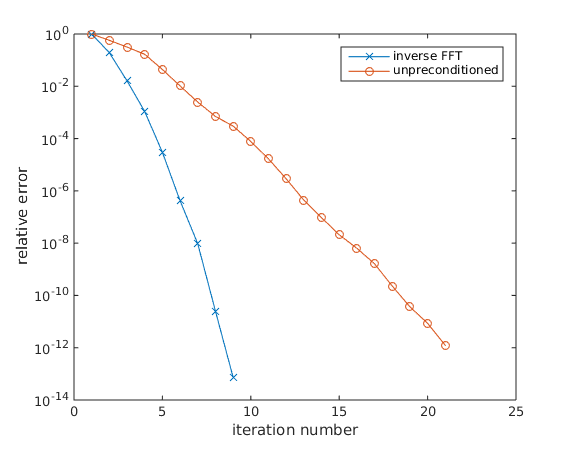}
    \includegraphics[height=200pt, width=220pt]{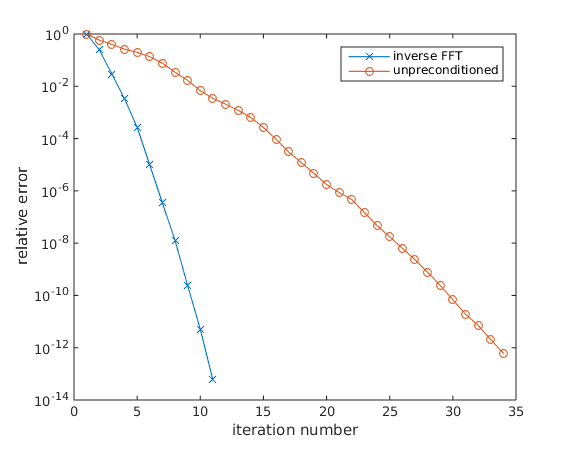}
    \includegraphics[height=200pt, width=220pt]{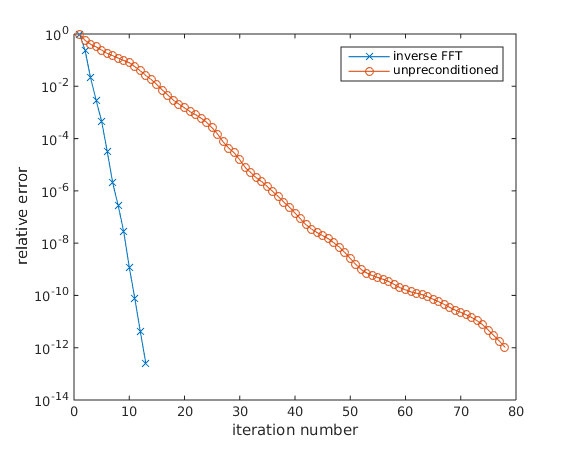}
    \includegraphics[height=200pt, width=220pt]{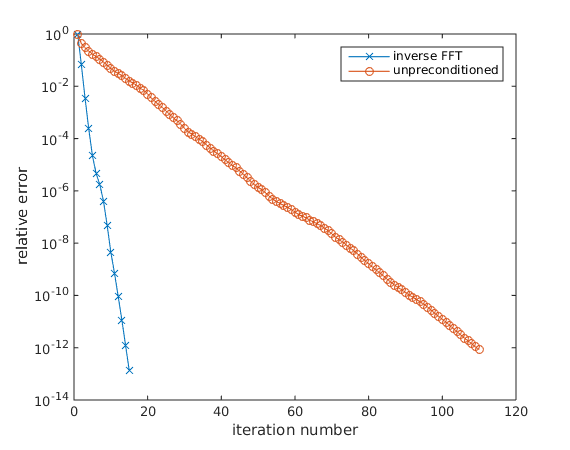}
    \caption{Convergence history for the diffusion preconditioned (blue $\times$) and unpreconditioned (orange $\circ$) \GMRES iterations with different scattering coefficients. Show from left and right are cases at $\mu_s=10.0$, $\mu_s = 20.0$, $\mu_s = 40.0$, and $\mu_s = 80.0$ respectively.}
    \label{FIG: EXP VI}
\end{figure}

\paragraph{Experiment VI.} We repeat the numerical simulations in Experiment V with the FFT-based preconditioner. The numerical results are summarized in Figure~\ref{FIG: EXP VI} and Table~\ref{TAB:EXP VI}. We observe very similar phenomena as we see in the case of the diffusion preconditioner, that is, the FFT-based is quite effective at reducing the total number of iterations it takes the algorithm to converge: the number of iterations to reach the convergence tolerence for the preconditioned \GMRES are respectively $9$ and $21$ ($\mu_s = 10.0$), $11$ and $34$ ($\mu_s=20.0$), $13$ and $78$ ($\mu_s=40.0$), and $15$ and $110$ ($\mu_s=80.0$). The compuational cost of each iteration for the preconditioned iteration is slightly higher than that of the unpreconditioned iteration. This results in significant gain in the overall computational efficiency of the preconditioned algorithm. In fact, the reduction of the overall cost is already happening at $\mu_s=10.0$ thanks for the dramatic reduction of the number of iterations to convergence. The main advantage of the FFT-based preconditioner over the diffusion-based preconditioner is that the number of preconditioned iteration to convergences grows very slow as the scattering coefficient increases, and the overall computational cost is much lower. %As a final remark, the above numerical preconditioner works when the total absorption coefficient varies little, not necessary to be a constant. If $\mu_t$ is given as a constant, then using FFT to compute the matrix-vector multiplication instead of low rank approximation might be faster for practical problem sizes.
\begin{table}[!htb]
     \centering
     \small
     \caption{The computational cost per iteration and the number of iterations to convergence for diffusion preconditioned and unpreconditioned \GMRES iterations under different scattering coefficients.}
     \vspace{0.2cm}
     \label{TAB:EXP VI}   
            \begin{tabular}{l*{10}{c}l}
                \hline
                $\mu_s$  & $\Tfg$ & $\Ifg$ & $\Tfpg$ & $\Ifpg$ \rule{0pt}{2.6ex}\rule[-1.2ex]{0pt}{0pt} \\
                \hline
                \rule{0pt}{4ex}$10.0$ &  $\ep{2.50}{+0}$ & $21$ & $\ep{2.91}{+0}$ & $9$ \\
                $20.0$ & $\ep{3.13}{+0}$  & $34$ & $\ep{3.54}{+0}$ & $11$ \\
                $40.0$ & $\ep{2.95}{+0}$  & $78$ & $\ep{3.29}{+0}$ & $13$\\
                $80.0$ & $\ep{3.11}{+0}$  & $110$ & $\ep{3.29}{+0}$ & $15$  %\\
                %$\eqref{EQ:VAR MUS}$ & $\ep{3.06}{+0}$  & $113$ & $\ep{3.27}{+0}$ & $17$ 
            \end{tabular}
\end{table}
%%%%%%%%%%%%%%%%%%%%%%%%%%%%%%%%%%%%%%%%%%%%%%%%%%%%%%%%%%%%%%%%%%
%%%%%%%%%%%%%%%%%%%%%%%%%%%%%%%%%%%%%%%%%%%%%%%%%%%%%%%%%%%%%%%%%%
\section{Concluding remarks}
\label{SEC:Concl}
%%%%%%%%%%%%%%%%%%%%%%%%%%%%%%%%%%%%%%%%%%%%%%%%%%%%%%%%%%%%%%%%%%
%%%%%%%%%%%%%%%%%%%%%%%%%%%%%%%%%%%%%%%%%%%%%%%%%%%%%%%%%%%%%%%%%%

To summarize, we presented in this work a fast numerical method for solving the equation of radiative transfer in isotropic media. The main idea of the method is to reformulate the ERT into an integral equation of the second type and then use the fast multipole technique to accelerate the solution of such an integral equation. Our numerical tests show that the algorithmic cost indeed scales linearly with respect to the size of the spatial component of the problem.

There are a few features of the method we proposed here. First, with the integral formulation, we avoid angular discretization of the ERT in the most expensive part of the solution process. This in principle allows us to handle large problems that would be hard to handle in, for instance, the discrete ordinate formulation, with limited RAM. Second, the kernel in our integral formulation of the ERT takes the same form for homogeneous and inhomogeneous media. Therefore, the algorithm we developed does not need to be modified going from homogeneous media problems to inhomogeneous media problems. This is quite different from existing fast multipole based methods. That said, in homogeneous media, the setup of our algorithm is relatively computationally inexpensive since the kernel in the corresponding integral equation is explicitly given. In inhomogeneous media, the setup requires the evaluation of the kernel for different $(\bx,\by)$ pairs that involves line integrals of the total absorption coefficients between $\bx$ and $\by$. This evaluation is more expensive than the homogeneous media case, but is still relatively low. In our implementation of the FMM algorithm, we cached all the calculations that involve the evaluation of the line integrals. This does not cause major storage problem since the number of Chebyshev interpolation nodes used in the implementation is always relatively small.

%In many practically relevant problems, we have coefficients that can be treated as periodic functions. Fast Fourier transform type of techniques can be used to accelerate the setup process of the algorithm. 

Let us also emphasize that, even though our formulation requires that the underlying medium to be isotropic, the internal and boundary source functions need not to be isotropic at all. In fact, the only thing that would have changed for the algorithm with an anisotropic source is the evaluation of $K(\mu_s^{-1}f)$.

We implemented a FFT-based and a diffusion-based preconditioning strategies for the solution of the linear integral equation involved in the calculation. Moreover, we observed from our numerical tests that the FMM approximation with a very small number of Chebyshev interpolation nodes already give relatively accuracy approximations to the true numerical solutions. This suggests that we can probably use the algorithm with small numbers of Chebyshev interpolation points as a preconditioning strategy for a general transport solver for more complicated problems. We are currently exploring in this direction.

To the best of our knowledge, what we proposed is the first algorithm for solving the ERT within the frame work of the fast multipole method. Our contribution is mainly on the introduction of the idea, not on the implementation of fast multipole methods. Indeed, our implementation is rather primitive which we believe can be improved, either by refining the current strategy or by exploring other approaches~\cite{YiBiZo-JCP03}. The study we have in this short paper is by no means enough to draw conclusions on every aspect of the algorithm, for instance how the algorithm benchmarks with existing methods. However, numerical simulations we have performed show that this is a promising method that is worth careful further investigated. We hope that this work can motivate more studies in this direction.

The generalization of our method to anisotropic media, that is, when the coefficients $\mu$ and $\mu_s$ depend on the angular variable $\bv$, is in general a quite challenging task. This has been done for some special forms of anisotropicity, for instance when the scattering kernel has only a small number of Fourier modes in the angular variable~\cite{FaAnYi-JCP19}. A recent study on the possibility of low-rank approximations to the integral kernel in anisotropic case can be found in~\cite{ReZhZh-arXiv19}.

%%%%%%%%%%%%%%%%%%%%%%%%%%%%%%%%%%%%%%%%%%%%%%%%%%%%%%%%%%%%%%%%%%
%%%%%%%%%%%%%%%%%%%%%%%%%%%%%%%%%%%%%%%%%%%%%%%%%%%%%%%%%%%%%%%%%%
\section*{Acknowledgments}
%%%%%%%%%%%%%%%%%%%%%%%%%%%%%%%%%%%%%%%%%%%%%%%%%%%%%%%%%%%%%%%%%%
%%%%%%%%%%%%%%%%%%%%%%%%%%%%%%%%%%%%%%%%%%%%%%%%%%%%%%%%%%%%%%%%%%

We would like to thank the anonymous referees for their constructive comments, including pointing out reference~\cite{GrHu-JCP02}, that help us improve the quality of this paper. During the revision of this paper, we were alerted about the paper~\cite{FaAnYi-JCP19} which generalized our algorithm to some simple anisotropic cases. This work is partially supported by the National Science Foundation through grant DMS-1620473.

%%%%%%%%%%%%%%%%%%%%%%%%%%%%%%%%%%%%%%%%%%%%%%%%%%%%%%%%%%%%%%%%%%
%%%%%%%%%%%%%%%%%%%%%%%%%%%%%%%%%%%%%%%%%%%%%%%%%%%%%%%%%%%%%%%%%%
\section*{Acknowledgments}
%%%%%%%%%%%%%%%%%%%%%%%%%%%%%%%%%%%%%%%%%%%%%%%%%%%%%%%%%%%%%%%%%%
%%%%%%%%%%%%%%%%%%%%%%%%%%%%%%%%%%%%%%%%%%%%%%%%%%%%%%%%%%%%%%%%%%

We would like to thank the anonymous referees for their constructive comments, including pointing out reference~\cite{GrHu-JCP02}, that help us improve the quality of this paper. During the revision of this paper, we were alerted about the paper~\cite{FaAnYi-JCP19} which generalized our algorithm to some simple anisotropic cases. This work is partially supported by the National Science Foundation through grant DMS-1620473.

%%%%%%%%%%%%%%%%%%%%%%%%%%%%%%%%%%%%%%%%%%%%%%%%%%%%%%%%%%%%%%%%%%
%%%%%%%%%%%%%%%%%%%%%%%%%%%%%%%%%%%%%%%%%%%%%%%%%%%%%%%%%%%%%%%%%%
%%%%%%%%%%%%%%%%%%%%%%%%%%%%%%%%%%%%%%%%%%%%%%%%%%%%%%%%%%%%%%%%%%
%%%%%%%%%%%%%%%%%%%%%%%%%%%%%%%%%%%%%%%%%%%%%%%%%%%%%%%%%%%%%%%%%%


{\small
\begin{thebibliography}{10}

\bibitem{AdLa-PNE02}
{\sc M.~L. Adams and E.~W. Larsen}, {\em Fast iterative methods for
  discrete-ordinates particle transport calculations}, Prog. Nucl. Energy, 40
  (2002), pp.~3--150.

\bibitem{AdNo-JCP98}
{\sc M.~L. Adams and P.~F. Nowak}, {\em Asymptotic analysis of a method for
  time- and frequency-dependent radiative transfer}, J. Comput. Physics, 146
  (1998), pp.~366--403.

\bibitem{Alcouffe-NSE77}
{\sc R.~E. Alcouffe}, {\em {Diffusion Synthetic Acceleration Methods for the
  Diamond-Differenced Discrete-Ordinates Equations}}, Nucl. Sci. Eng., 64
  (1977), p.~344.

\bibitem{AnLa-JCP01}
{\sc D.~Y. Anistratov and E.~W. Larsen}, {\em Nonlinear and linear
  $\alpha$-weighted methods for particle transport problems}, J. Comput. Phys.,
  173 (2001), pp.~664--684.

\bibitem{Arridge-IP99}
{\sc S.~R. Arridge}, {\em Optical tomography in medical imaging}, Inverse
  Probl., 15 (1999), pp.~R41--R93.

\bibitem{Asadzadeh-SIAM98}
{\sc M.~Asadzadeh}, {\em A finite element method for the neutron transport
  equation in an infinite cylindrical domain}, SIAM J. Numer. Anal., 35 (1998),
  pp.~1299--1314.

\bibitem{BaCaLiRe-IP07}
{\sc G.~Bal, L.~Carin, D.~Liu, and K.~Ren}, {\em Experimental validation of a
  transport-based imaging method in highly scattering environments}, Inverse
  Problems, 23 (2007), pp.~2527--2539.

\bibitem{BaRe-IP05}
{\sc G.~Bal and K.~Ren}, {\em Atmospheric concentration profile reconstructions
  from radiation measurements}, Inverse Problems, 21 (2005), pp.~153--168.

\bibitem{BaRe-SIAM08}
\leavevmode\vrule height 2pt depth -1.6pt width 23pt, {\em Transport-based
  imaging in random media}, SIAM J. Appl. Math., 68 (2008), pp.~1738--1762.

\bibitem{BeGr-WMMEP97}
{\sc R.~Beatson and L.~Greengard}, {\em A short course on fast multipole
  methods}, in Wavelets, Multilevel Methods and Elliptic PDEs, Oxford
  University Press, 1997, pp.~1--37.

\bibitem{BhSp-JCP07}
{\sc K.~Bhan and J.~Spanier}, {\em Condensed history {Monte Carlo} methods for
  photon transport problems}, J. Comput. Phys., 225 (2007), pp.~1673--1694.

\bibitem{Boman-Thesis07}
{\sc E.~Boman}, {\em Radiotherapy Forward and Inverse Problem Applying
  Boltzmann Transport Equations}, PhD thesis, University of Kuopio, Filand,
  Kuopio, Filand, 2007.

\bibitem{BoGa-IP16}
{\sc L.~Borcea and J.~Garnier}, {\em Derivation of a one-way radiative transfer
  equation in random media}, Phys. Rev. E, 93 (2016).
\newblock 022115.

\bibitem{BoLaAd-JCP92}
{\sc C.~Borgers, E.~W. Larsen, and M.~L. Adams}, {\em {The asymptotic diffusion
  limit of a linear discontinuous discretization of a 2-dimensional linear
  transport equation}}, J. Comp. Phys., 98(2) (1992), pp.~285--300.

\bibitem{BrGi-JCP12}
{\sc J.~Bremer and Z.~Gimbutas}, {\em A {Nystr{\"o}m} method for weakly
  singular integral operators on surfaces}, J. Comput. Phys., 231 (2012),
  pp.~4885--4903.

\bibitem{BrMoRa-JCP14}
{\sc D.~E. Bruss, J.~E. Morel, and J.~C. Ragusa}, {\em {$S_2SA$}
  preconditioning for the {$S_N$} equations with strictly nonnegative spatial
  discretization}, J. Comput. Phys., 273 (2014), pp.~706--719.

\bibitem{CeBaBeAi-TTSP99}
{\sc C.~Cecchi-Pestellini, L.~Barletti, A.~Belleni-Morante, and S.~Aiello},
  {\em Radiative transfer in the stochastic interstellar medium}, Trans. Theor.
  Stat. Phys., 28 (1999), pp.~199--228.

\bibitem{ChGrRo-JCP99}
{\sc H.~Cheng, L.~Greengard, and V.~Rokhlin}, {\em A fast adaptive multipole
  algorithm in three dimensions}, J. Comp. Phys., 155 (1999), pp.~468--498.

\bibitem{ChJiMi-Book01}
{\sc W.~C. Chew, J.~M. Jin, and E.~Michielssen}, eds., {\em Fast and Efficient
  Algorithms in Computational Electromagnetics}, Artech House Publishers, 2001.

\bibitem{DaLi-Book93-6}
{\sc R.~Dautray and J.-L. Lions}, {\em Mathematical {Analysis} and {Numerical}
  {Methods} for {Science} and {Technology}, {Vol VI}}, Springer-Verlag, Berlin,
  1993.

\bibitem{DeVo-JCP02}
{\sc A.~Dedner and P.~Vollm\"oller}, {\em An adaptive higher order method for
  solving the radiation transport equation on unstructured grids}, J. Comput.
  Phys., 178 (2002), pp.~263--289.

\bibitem{DeThUr-JCP12}
{\sc J.~D. Densmore, K.~G. Thompson, and T.~J. Urbatsch}, {\em A hybrid
  transport-diffusion {Monte Carlo} method for frequency-dependent
  radiative-transfer simulations}, J. Comput. Phys., 231 (2012),
  pp.~6924--6934.

\bibitem{DiRe-JCP14}
{\sc T.~Ding and K.~Ren}, {\em Inverse transport calculations in optical
  imaging with subspace optimization algorithms}, J. Comput. Phys., 273 (2014),
  pp.~212--226.

\bibitem{DuKl-JCP02}
{\sc B.~Dubroca and A.~Klar}, {\em Half-moment closure for radiative transfer
  equations}, J. Comput. Phys., 180 (2002), pp.~584--596.

\bibitem{DuGuRo-SIAM96}
{\sc A.~Dutt, M.~Gu, and V.~Rokhlin}, {\em Fast algorithms for polynomial
  interpolation, integration, and differentiation}, SIAM J. Numer. Anal., 33
  (1996), pp.~1689--1711.

\bibitem{Edstrom-SIAM05}
{\sc P.~Edstr\"{o}m}, {\em A fast and stable solution method for the radiative
  transfer problem}, SIAM Rev., 47 (2005), pp.~447--468.

\bibitem{FaAnYi-JCP19}
{\sc Y.~Fan, J.~An, and L.~Ying}, {\em Fast algorithms for integral
  formulations of steady-state radiative transfer equation}, J. Comput. Phys.,
  380 (2019), pp.~191--211.

\bibitem{FoDa-JCP09}
{\sc W.~Fong and E.~Darve}, {\em The black-box fast multipole method}, J.
  Comput. Phys., 228 (2009), pp.~8712--8725.

\bibitem{FrKlLaYa-JCP07}
{\sc M.~Frank, A.~Klar, E.~W. Larsen, and S.~Yasuda}, {\em Time-dependent
  simplified {$P_N$} approximation to the equations of radiative transfer}, J.
  Comput. Phys., 226 (2007), pp.~2289--2305.

\bibitem{GaZh-TTSP09}
{\sc H.~Gao and H.~Zhao}, {\em A fast forward solver of radiative transfer
  equation}, Trans. Theor. Stat. Phys., 38 (2009), pp.~149--192.

\bibitem{GoLi-JCP12}
{\sc W.~F. Godoy and X.~Liu}, {\em Parallel {Jacobian-free Newton Krylov}
  solution of the discrete ordinates method with flux limiters for {3D}
  radiative transfer}, J. Comput. Phys., 231 (2012), pp.~4257--4278.

\bibitem{Graziani-Book06}
{\sc F.~Graziani}, ed., {\em Computational Methods in Transport}, Springer,
  2006.

\bibitem{GrHu-JCP02}
{\sc L.~Greengard and J.~Huang}, {\em A new version of the fast multipole
  method for screened {Coulomb} interactions in three dimensions}, J. Comput.
  Phys., 180 (2002), pp.~642--658.

\bibitem{GrRo-JCP87}
{\sc L.~Greengard and V.~Rokhlin}, {\em A fast algorithm for particle
  simulations}, J. Comput. Phys., 73 (1987), pp.~325--348.

\bibitem{GrRo-AN97}
\leavevmode\vrule height 2pt depth -1.6pt width 23pt, {\em A new version of the
  fast multipole method for the {Laplace} equation in three dimensions}, Acta
  Numer., 6 (1997), pp.~229--270.

\bibitem{GrSc-JCP11}
{\sc K.~Grella and C.~Schwab}, {\em Sparse tensor spherical harmonics
  approximation in radiative transfer}, J. Comput. Phys., 230 (2011),
  pp.~8452--8473.

\bibitem{GuKa-SIAM10}
{\sc J.-L. Guermond and G.~Kanschat}, {\em Asymptotic analysis of upwind {DG}
  approximation of the radiative transport equation in the diffusive limit},
  SIAM J. Numer. Anal., 48 (2010), pp.~53--78.

\bibitem{HaSpVe-SIAM07}
{\sc C.~K. Hayakawa, J.~Spanier, and V.~Venugopalan}, {\em Coupled
  forward-adjoint monte carlo simulations of radiative transport for the study
  of optical probe design in heterogeneous tissues}, SIAM J. Appl. Math., 68
  (2007), pp.~253--270.

\bibitem{HeGr-AJ41}
{\sc L.~G. Henyey and J.~L. Greenstein}, {\em Diffuse radiation in the galaxy},
  Astrophys. J., 90 (1941), pp.~70--83.

\bibitem{Hermeline-JCP16}
{\sc F.~Hermeline}, {\em A discretization of the multigroup {$P_N$} radiative
  transfer equation on general meshes}, J. Comput. Phys., 313 (2016),
  pp.~549--582.

\bibitem{Hochstadt-Book89}
{\sc H.~Hochstadt}, {\em Integral Equations}, Wiley, 1989.

\bibitem{JiScPaJi-PMB12}
{\sc X.~Jia, J.~Sch\"umann, H.~Paganetti, and S.~B. Jiang}, {\em {GPU}-based
  fast monte carlo dose calculation for proton therapy}, Phys. Med. Biol., 57
  (2012), pp.~7783--7797.

\bibitem{JiPaTo-SIAM00}
{\sc S.~Jin, L.~Pareschi, and G.~Toscani}, {\em Uniformly accurate diffusive
  relaxation schemes for multiscale transport equations}, SIAM J. Numer. Anal.,
  38 (200), pp.~913--936.

\bibitem{KaRa-SIAM14}
{\sc G.~Kanschat and J.~C. Ragusa}, {\em A robust multigrid preconditioner for
  {SNDG} approximation of monochromatic, isotropic radiation transport
  problems}, SIAM J. Sci. Comput., 36 (2014), pp.~2326--2345.

\bibitem{KiMo-SIAM02}
{\sc A.~D. Kim and M.~Moscoso}, {\em Chebyshev spectral methods for radiative
  transfer}, SIAM J. Sci. Comput., 23 (2002), pp.~2075--2095.

\bibitem{KiMo-IP06}
\leavevmode\vrule height 2pt depth -1.6pt width 23pt, {\em Radiative transport
  theory for optical molecular imaging}, Inverse Problems, 22 (2006),
  pp.~23--42.

\bibitem{KiBeGoMo-JCP10}
{\sc M.~Kindelan, F.~Bernal, P.~Gonz\'alez-Rodr\'iguez, and M.~Moscoso}, {\em
  Application of the {RBF} meshless method to the solution of the radiative
  transport equation}, J. Comput. Phys., 229 (2010), pp.~1897--1908.

\bibitem{Kress-Book99}
{\sc R.~Kress}, {\em Linear Integral Equations}, Applied Mathematical Sciences,
  Springer-Verlag, New York, 2nd~ed., 1999.

\bibitem{LaMcHa-JCP16}
{\sc V.~M. Laboure, R.~G. McClarren, and C.~D. Hauck}, {\em Implicit filtered
  {$P_N$} for high-energy density thermal radiation transport using
  discontinuous {Galerkin} finite elements}, J. Comput. Phys., 321 (2016),
  pp.~624--643.

\bibitem{LaPaSe-Book03}
{\sc B.~Lapeyre, E.~Pardoux, and R.~Sentis}, {\em Introduction to Monte-Carlo
  Methods for Transport and Diffusion Equations}, Oxford University Press,
  2003.

\bibitem{Larsen-TTSP84}
{\sc E.~W. Larsen}, {\em Diffusion-synthetic acceleration methods for
  discrete-ordinates problems}, Transport Theory and Statistical Physics, 13
  (1984), pp.~107--126.

\bibitem{Larsen-NSE88}
{\sc E.~W. Larsen}, {\em Neutronics methods for thermal radiative transfer},
  Nucl. Sci. Eng., 100 (1988), pp.~255--259.

\bibitem{Larsen-TTSP88}
\leavevmode\vrule height 2pt depth -1.6pt width 23pt, {\em Solution of
  three-dimensional inverse transport problems}, Trans. Theor. Stat. Phys., 17
  (1988), pp.~147--167.

\bibitem{LaThKlSeGo-JCP02}
{\sc E.~W. Larsen, G.~Th\"ommes, A.~Klar, M.~Sea\"idd, and T.~G\"otze}, {\em
  Simplified {PN} approximations to the equations of radiative heat transfer
  and applications}, J. Comput. Phys., 183 (2002), pp.~652--675.

\bibitem{Lesaint-FEP86}
{\sc P.~Lesaint}, {\em Finite element methods for the transport equation}, in
  Finite Elements in Physics, North-Holland, Amsterdam, 1987.

\bibitem{LeMi-Book93}
{\sc E.~E. Lewis and W.~F. Miller}, {\em Computational Methods of Neutron
  Transport}, American Nuclear Society, La Grange Park, IL, 1993.

\bibitem{MaReSt-SIAM00}
{\sc T.~A. Manteuffel, K.~Ressel, and G.~Starke}, {\em A boundary functional
  for the least-squares finite-element solution of the neutron transport
  equation}, SIAM J. Numer. Anal., 37(2) (2000), pp.~556--586.

\bibitem{MaRo-SIAM07}
{\sc P.~Martinsson and V.~Rokhlin}, {\em An accelerated kernel-independent fast
  multipole method in one dimension}, SIAM J. Sci. Comput., 29 (2007),
  pp.~1160--1178.

\bibitem{Mokhtar-Book97}
{\sc M.~Mokhtar-Kharroubi}, ed., {\em {Mathematical} {Topics} in {Neutron}
  {Transport} {Theory}: {New} {Aspects}}, World Scientific, Singapore, 1997.

\bibitem{MoRaAdKa-TTSP13}
{\sc J.~E. Morel, J.~C. Ragusa, M.~L. Adams, and G.~Kanschat}, {\em Asymptotic
  {PN}-equivalent {SN}+1 equations}, Transport Theory Stat. Phys., 42 (2013),
  pp.~3--20.

\bibitem{OlDe-PNE98}
{\sc S.~Oliveira and Y.~Deng}, {\em Preconditioned {Krylov} subspace methods
  for transport equations}, Progress in Nuclear Energy, 33 (1998),
  pp.~155--174.

\bibitem{PaHo-ANE02}
{\sc B.~W. Patton and J.~P. Holloway}, {\em Application of preconditioned
  {GMRES} to the numerical solution of the neutron transport equation}, Annals
  of Nuclear Energy, 29 (2002), pp.~109--136.

\bibitem{RaGuKa-JCP12}
{\sc J.~Ragusa, J.-L. Guermond, and G.~Kanschat}, {\em A robust
  {Sn}-{DG}-approximation for radiation transport in optically thick and
  diffusive regimes}, J. Comput. Phys., 231 (2012), pp.~1947--1962.

\bibitem{Ren-CiCP10}
{\sc K.~Ren}, {\em Recent developments in numerical techniques for
  transport-based medical imaging methods}, Commun. Comput. Phys., 8 (2010),
  pp.~1--50.

\bibitem{ReAbBaHi-OL04}
{\sc K.~Ren, G.~S. Abdoulaev, G.~Bal, and A.~H. Hielscher}, {\em Algorithm for
  solving the equation of radiative transfer in the frequency domain}, Optics
  Lett., 29 (2004), pp.~578--580.

\bibitem{ReBaHi-SIAM06}
{\sc K.~Ren, G.~Bal, and A.~H. Hielscher}, {\em Frequency domain optical
  tomography based on the equation of radiative transfer}, SIAM J. Sci.
  Comput., 28 (2006), pp.~1463--1489.

\bibitem{ReZhZh-arXiv19}
{\sc K.~Ren, H.~Zhao, and Y.~Zhong}, {\em Separability of the kernel function
  in an integral formulation for anisotropic radiative transfer equation},
  Submitted,  (2019).

\bibitem{RoViBo-PNE11}
{\sc B.~D. Rodriguez, M.~T. Vilhena, and B.~E.~J. Bodmann}, {\em An overview of
  the {Boltzmann} transport equation solution for neutrons, photons and
  electrons in {Cartesian} geometry}, Progress in Nuclear Energy, 53 (2011),
  pp.~1119--1125.

\bibitem{RyPaKe-WM96}
{\sc L.~Ryzhik, G.~Papanicolaou, and J.~B. Keller}, {\em Transport equations
  for elastic and other waves in random media}, Wave Motion, 24 (1996),
  pp.~327--370.

\bibitem{Saad-Book03}
{\sc Y.~Saad}, {\em Iterative Methods for Sparse Linear Systems}, SIAM,
  Philadelphia, 2nd~ed., 2003.

\bibitem{SpKuCh-JQSRT01}
{\sc R.~J.~D. Spurr, T.~P. Kurosu, and K.~V. Chance}, {\em A linearized
  discrete ordinate radiative transfer model for atmospheric remote-sensing
  retrieval}, J. Quant. Spectrosc. Radiat. Transfer, 68 (2001), pp.~689--735.

\bibitem{TiHeToSiAlPyUl-PMB08}
{\sc L.~Tillikainen, H.~Helminen, T.~Torsti, S.~Siljam\"aki, J.~Alakuijala,
  J.~Pyyry, and W.~Ulmer}, {\em A 3d pencil-beam-based superposition algorithm
  for photon dose calculation in heterogeneous media}, Phys. Med. Biol., 53
  (2008), pp.~3821--3839.

\bibitem{TuFrDuKl-JCP04}
{\sc R.~Turpault, M.~Frank, B.~Dubroca, and A.~Klar}, {\em Multigroup half
  space moment approximations to the radiative heat transfer equations}, J.
  Comput. Phys., 198 (2004), pp.~363--371.

\bibitem{UeLa-JCP98}
{\sc T.~Ueki and E.~Larsen}, {\em {A kinetic theory for nonanalog {Monte Carlo}
  particle transport algorithms: exponential transform with angular biasing in
  planar-geometry anisotropically scattering media}}, J. Comp. Phys., 145
  (1998), pp.~406--431.

\bibitem{Vainikko-Book93}
{\sc G.~Vainikko}, {\em {Multidimensional Weakly Singular Integral Equations}},
  Springer-Verlag, 1993.

\bibitem{WaUe-ASS89}
{\sc A.~P. Wang and S.~Ueno}, {\em An inverse problem in a three-dimensional
  radiative transfer}, Astrophys. Space Sci., 155 (1989), pp.~105--111.

\bibitem{YiBiZo-JCP03}
{\sc L.~Ying, G.~Biros, and D.~Zorin}, {\em A kernel-independent adaptive fast
  multipole method in two and three dimensions}, J. Comput. Phys., 196 (2003),
  pp.~591--626.

\end{thebibliography}
}
\end{document}